\documentclass[journal,twoside]{IEEEtran}
\ifCLASSINFOpdf
\else
\fi
\usepackage{amsmath}
\usepackage{array}
\usepackage{url}

\usepackage{tabularx}
\usepackage{mathrsfs}
\usepackage{bm}
\usepackage{siunitx}
\usepackage{todonotes}
\usepackage{amssymb}
\usepackage{dsfont}
\usepackage{bbold}
\usepackage{cite}
\usepackage{enumerate}

\newtheorem{proposition}{Proposition}

\newtheorem{remark}{Remark}
\newtheorem{assumption}{Assumption}

\newcommand{\nullvec}[1]{\mathbb 0_{#1}}

\newenvironment{neu}{}{}

\begin{document}
%
\title{Optimal Distributed Frequency and Voltage Control for Zonal Electricity Markets}
%
%
%
\author{Lukas~Kölsch,
        Lena~Zellmann,
        Rishabh~Vyas,
        Martin~Pfeifer,
        and~Sören~Hohmann,~\IEEEmembership{Member,~IEEE}
\thanks{L. Kölsch, L. Zellmann, R. Vyas, M. Pfeifer, and S. Hohmann are with the Department of Electrical Engineering and Information Technology, Karlsruhe Institute of Technology, Karlsruhe, Germany. (e-mail: lukas.koelsch@kit.edu).}\thanks{This work was funded by the Deutsche Forschungsgemeinschaft (DFG, German Research Foundation) -- project number 360464149.}}

%
%

\markboth{}%
{}
%



\maketitle

\begin{abstract}
Zonal pricing is a well-suited mechanism to incentivize grid-supporting behavior of profit-maximizing producers and consumers operating on a large-scale power system.
In zonal electricity markets, local system operators create individual price zones, 
which provide appropriate price signals depending on local grid conditions such as an excess or shortage of electrical energy in certain regions.
In this paper, a \begin{neu} real-time\end{neu} zonal pricing controller for AC power networks is presented that ensures
frequency and voltage stability as well as Pareto efficiency of the resulting closed-loop equilibria.
Based on a dynamic network model \begin{neu}that takes\end{neu} line losses and power exchange with adjacent price zones \begin{neu}into account\end{neu}, distributed \begin{neu}{continuous-time}\end{neu} control laws are derived which require only neighbor-to-neighbor communication.
\begin{neu}Application to a real-time congestion management strategy
	illustrates how spatially and temporally differentiated prices enable grid-supportive operation
	of the participants without interventions by a superordinate control authority.
	Effectiveness of different zonal pricing concepts compared to an isolated grid operation is demonstrated
	through simulations on the IEEE-57 bus system.\end{neu}
\end{abstract}


%
\IEEEpeerreviewmaketitle

\begin{neu}\section*{Nomenclature}
\subsection*{Abbreviations}
\noindent\begin{tabularx}{0.89\columnwidth}{@{}p{1.1cm}X}
	CC & Cell Coordinator\\
	KKT & Karush-Kuhn-Tucker\\
	MU & Monetary Unit\\
	OPF & Optimal Power Flow\\
	PPO & Power Plant Operator\\
	RES & Renewable Energy Source\\
	SG & Synchronous Generator\\
	VPP & Virtual Power Plant\\
	WoC & Web-of-Cells
	\end{tabularx}
\subsection*{Sets and Tuples}
\noindent\begin{tabularx}{0.89\columnwidth}{@{}p{1.1cm}X}
		$\mathcal E_c$ & set of communication links\\
	$\mathcal E_p, \widehat{\mathcal E}_p$ & set of lines/ inter-cell lines\\
	$\mathscr G_c$ & communication graph \\
	${\mathscr G}_p$ & graph of physical system\\
	$\mathscr G_z$ & graph of cell interconnection\\
	$\mathcal N$ & set of neighboring nodes\\
	$\mathcal P$ & set of PPOs\\
	$\mathcal V$ & set of buses \\
	$\mathcal Z$ & set of cells
\end{tabularx}
\subsection*{Parameters}
\noindent\begin{tabularx}{0.89\columnwidth}{@{}p{1.1cm}X}
	$A_i$ & damping coefficient \\
	$\bm B$ & negative susceptance matrix \\
	$\mathcal B$ & Laplacian matrix\\
	$b_i$ & shunt susceptance\\
	$\bm D$ & node-edge incidence matrix \\
	$\bm G$ & negative conductance matrix \\
	$g_i$ & shunt conductance\\
	$L_i$ & angular momentum deviation\\
	$M_i$ & moment of inertia\\
	$X_{d,i}, X_{d,i}'$ & d-axis synchronous/ transient reactance \\
	$\tau_i$ & time constant 
\end{tabularx}
\subsection*{Variables}
\noindent\begin{tabularx}{0.89\columnwidth}{@{}p{1.1cm}X}
	$\mathtt C_{i}$ & cost \\
	$\mathcal C_m$ & congestion rate\\
	$\mathscr L$ & Lagrangian\\
	$\mathtt P_i$ & profit\\
	$P_{ij}$ & sending-end active power flow \\
	$p_i,q_i$ & active/ reactive power injection \\
	$\mathtt U_{i}$ & utility \\
	$U_i$ & magnitude of transient internal voltage\\
	$U_{f,i}$ & magnitude of excitation voltage\\
	$\kappa_k$ & participation factor of cell-specific prices \\
	$\Lambda_k$ & price in cell $k$\\ 
	$\lambda_i$ &price at node $i$\\
	$\mu,\nu$ & Lagrange multiplier\\
	$\vartheta_{ij}$ & bus voltage angle difference\\
	$\Phi$ & transmission loss\\
	$\omega_i$ & frequency deviation from nominal value
\end{tabularx}
\subsection*{Indices and Exponents}
\noindent\begin{tabularx}{0.89\columnwidth}{@{}p{1.1cm}X}
	$\square_{\mathcal G}$ & SG node\\
	$\square_{g}$ & generation \\
	$\square_{\mathcal I}$ & inverter node\\
	$\square_{k}$ & cell \\
	$\square_{\mathcal L}$ & load node\\
	$\square_{\ell}$ & consumption \\
	$\square_{\pi}$ & PPO\\
	$\square^{\star},\square^{\sharp}$ & optimizer\\
	$\overline{\square},\underline{\square}$ & upper limit/ lower limit
\end{tabularx}\end{neu}

\section{Introduction}
\IEEEPARstart{T}{he} \begin{neu} vast majority of today's power system is governed by an interplay of profit-maximizing prosumers and regulated system operators.
The reason for this coexistence is the principle of unbundling, according to which the responsibilities for energy supply and grid operation must be economically separated \cite{EU.2009b}.\end{neu}
As a result of the worldwide trend towards more renewable power generation and the displacement of large conventional power plants, there is an increasing number of small-scale generation.
A major challenge of this trend is the growing complexity and heterogeneous nature of the future power system consequent from different types of distributed generation units, such as conventional synchronous generators (SGs), renewable energy sources (RESs), and power consuming loads \cite{Kotsampopoulos.2013}.
\begin{neu}
Consequently, the number of competitive grid participants is increasing, which in turn hampers a centralized monitoring and coordination by a superordinate instance.
\end{neu}
From an optimization-theoretic perspective, the interaction between competitive grid participants can be interpreted as a \begin{neu} network-constrained non-cooperative game\end{neu}, where the objective functions are given by the individual profits and the constraints are imposed by physical laws as well as technical and/ or operational regulations.

\begin{neu}
	In such a competitive environment, \emph{real-time pricing} is considered to have become an important component to enforce a grid-supportive behavior of selfishly motivated players \cite{Li.2015,Wang.2015b,Dutta.2017}.\todo[disable]{Wang2015b kann evtl. ganz rausfliegen.}
	It allows to reflect the external costs caused by the network, e.g. frequency deviations due to global demand-supply mismatch or transmission line congestion in real-time prices
	and thus provides short-term incentives to \begin{neu}elastic\end{neu} producers and consumers for restoring the demand-supply balance or clearing the congestion without the need of after-market interventions. 
	A major challenge for the deployment of a pricing mechanism is how to
	appropriately internalize the external costs without having to rely on a centralized optimization problem
	{while satisfying} the privacy constraints of participants.
	
\end{neu}

Yet, the existing \begin{neu}price-based control frameworks\end{neu} mainly deal with centralized or sometimes distributed optimization problems, but the majority of work considers only one objective function.
For this problem class, the primal-dual gradient method \cite{Feijer.2010} allows to derive continuous-time feedback controllers, which has
%
already been successfully applied in particular aspects of power system control such as frequency regulation \cite{Stegink.2017,Mallada.2017,Dorfler.2017}, voltage regulation \cite{Magnusson.2017} or optimal power flow \cite{Ma.2013,Zhang.2015}. A detailed discussion of current research work on optimization-based control of power systems in continuous time can be found in \cite{Dorfler.2019}. 
By applying the primal-dual gradient method, the Lagrange multipliers of the local power balance constraints are associated with the marginal costs of generation, which turns the gradient ascent of the dual variables into a \begin{neu}continuous-time\end{neu} 
dynamic pricing controller (cf. \cite{Liu.2009c,Molzahn.2017}).

In the context of \begin{neu}non-cooperative frameworks
	for frequency control\end{neu}, the literature is scarce. 
Two recent papers \cite{Stegink.2019b} and \cite{Cherukuri.2020b}
develop a Bertrand competition model between price-setting generators
resulting in a continuous-time bidding process against a centralized system operator, which is shown to provide economic efficiency.
The authors of \cite{Persis.2019} propose a distributed dynamic pricing mechanism by feedback optimization
for a Cournot model of competition between price-setting generators and price-taking (i.e. \begin{neu}elastic\end{neu}) loads.
However, both approaches rely on the simplifying assumptions of constant voltage magnitudes and zero line losses, i.e. they do not account for the physical system.
\begin{neu}As a consequence, zero frequency deviation can be accomplished only if all line resistances of the system are zero
and furthermore, power flow constraints can only be included if the topology of the network is radial.
\end{neu}

\begin{neu}
However, apart from global supply-demand balancing, regionally uniform prices are in general unable to reflect \emph{local} scarcity 
since geo-spatial network effects such as power flow limits cannot be internalized \cite{Trepper.2015}.
In particular, the methods discussed so far do not prevent from after-market redispatch by the system operator due to transmission congestion. 
As a remedy, \emph{zonal pricing} is considered to reduce structural and regional mismatches in generation and consumption by signalling local network conditions and thus providing
short- and long-term investment signals to deploy generation capacity in load pockets \cite{Staudt.2019}.
Besides numerous research contributions on a rather strategic level which typically analyze the long-term perspective, a few papers deal with real-time zonal pricing schemes to address congestion management.
{The authors of} \cite{MahmoudianEsfahani.2016} and \cite{MahmoudianEsfahani.2017} present centralized strategies for cost-efficient generator rescheduling based on an offline calculation of sensitivity factors.
In \cite{Shiltz.2016,Shiltz.2019}, active power flow limits are included in a centralized real-time economic controller, while \cite{Zhang.2015,Zhang.2015a} present a fully distributed control scheme. 
However, all of the above pricing approaches rely on simplified network models such as DC-OPF or decoupled AC-OPF, where the crucial effects of transmission losses are neglected at market clearing, thus again provoking steady-state frequency errors.
\end{neu}

\begin{neu}
	An intriguing approach 
towards\end{neu} a fully decentralized power system \begin{neu}was proposed with\end{neu} the \emph{web of cells} (WoC) concept \cite{Martini.2017,Lehmann.2019b}.
\begin{neu} In this concept,\end{neu} the overall system is divided into interconnected subsystems, which are each monitored and controlled by a techno-economic \emph{cell coordinator} (CC). 
\begin{neu}By activating and maintaining automatic control mechanisms,
CCs
are supposed to solve local problems locally \cite{Rikos.2017,Cabiati.2018}.\end{neu}
\begin{neu} In particular, each CC\end{neu}
is responsible for providing frequency and voltage stability within its own cell by applying zonal (i.e. cell-specific) price incentives \begin{neu}\cite{MacDougall.2017}\end{neu}.
Due to their connection via transmission lines, the cells \begin{neu}{physically}\end{neu} interact with each other.
However, to the best of the authors' knowledge, there exists no explicit price-based feedback control strategy \begin{neu}complying with the\end{neu} WoC \begin{neu}concept\end{neu} to date.
\paragraph*{Statement of Contributions}
In this paper, we derive a distributed, continuous-time zonal pricing controller for lossy AC power networks.
For this purpose, we combine a feedback controller for \begin{neu}real-time\end{neu} optimization with the WoC concept, enabling incentive-based regional price differences as well as power flows across cell boundaries.
In line with the WoC concept, our work assumes a mixture of price-taking power plant operators (PPOs), price-setting CCs, and in\begin{neu}elastic\end{neu} consumers. 
The simultaneous execution of distributed controllers results in a 
dynamic pricing procedure, which ensures both Pareto-efficient allocations as well as frequency and voltage stability.
\begin{neu}
	{While the absolute values of the resulting prices depend on the current supply-demand conditions and transmission losses,
	the desired ratios of the individual cell-specific prices can be chosen freely within our approach.}
	In the second part of the paper, we exploit this degree of freedom by deriving an exemplary feedback control strategy for flow-based congestion management, which automatically adjusts the local price differences in the event of heavily loaded inter-cell lines in certain regions.	
	In contrast to existing approaches
	based on offline calculation of generator shift keys and power transmission distribution factors, the presented control strategy does not require a global snapshot of the grid state. Instead, all relevant parameters are calculated online based on local measurements and neighbor-to-neighbor communication between adjacent CCs.\end{neu}
Since the individual resources of a PPO may be located at geographically distant nodes in different cells, \begin{neu}our approach\end{neu} allows the integration of Virtual Power Plant (VPP) operators aggregating a large number of distributed small-scale RESs.
The cell-based control architecture \begin{neu}enables\end{neu} a global exchange of power, while measurement and control information is only shared locally.
In particular, none of the individual profits, costs, \begin{neu}or private constraints\end{neu} need to be disclosed by the network participants. 
The controller design is applicable for power systems with \begin{neu}a mixture of\end{neu} conventional and renewable \begin{neu}{power generation resources}\end{neu}.


\todo[disable]{noch irgendwo erwähnen:[wir verfolgen das konzept von locational marginal prices, d.h.: Die Preise sind flexibel und konvergieren anschließend zu einem gemeinsamen Wert. Je nach Konfiguration kann das Verhältnis der zellspezifischen Preise unbestimmt sein, oder aber zur Laufzeit ein bestimmtes Verhältnis annehmen] }

\paragraph*{Paper Organization}
Section \ref{ch:woc-physical-model} presents the \begin{neu}physical\end{neu} network model. 
Section \ref{ch:woc-problemformulation} formulates the optimization problems of the individual network participants and provides a deduction of continuous-time distributed controllers.
Section \ref{ch:woc-multiple-cells} extends the controller towards a WoC scenario with mutually related zonal prices and studies the Pareto efficiency of the closed-loop equilibrium. 
\begin{neu}In section \ref{ch:application}, automatic regulation of the zonal price differences
 is employed to develop a distributed controller for real-time congestion management.\end{neu}
A simulation study using the IEEE 57-bus system is carried out in section \begin{neu}\ref{ch:woc-simulation}\end{neu}.
Section \begin{neu}\ref{ch:woc-conclusion}\end{neu} gives the conclusion and proposes directions for future research.
\paragraph*{Notation}
Vectors and matrices are written in boldface.
All vectors defined in the paper are column vectors $\bm a = \mathrm{col}_i\{a_i\}=\mathrm{col}\{a_1,a_2,\ldots\}$ with elements $a_i$, $i=1,2,\ldots$.
All-zeros and all-ones vectors with $n$ entries are denoted by $\nullvec{n}$ and $\mathds 1_n$, respectively.
The $(n \times n)$-identity matrix is denoted by $\bm I_n$.
Positive definiteness of a matrix is denoted by $\succ 0$.
For $\bm a \in \mathds R^n$, we write $\bm a > \nullvec{n}$ if each component in $\bm a$ is greater than zero.
Upper and lower bounds are denoted by $\overline{\square}$ and $\underline{\square}$, respectively.
For a given $\mu \geq 0$, define
\begin{align}
\left\langle x\right\rangle_\mu^+ := \left\{ \begin{array}{ll} x, &\quad \mu > 0 \vee x \geq 0, \\ 0, &\quad \text{\begin{neu}otherwise\end{neu}.} \end{array}\right. \label{plusoperator}
\end{align} 
If $\bm x$ and $\bm \mu$ are vectors of the same size $i \in \mathds N$, then \eqref{plusoperator} can be applied component-wise, i.e.
$\langle \bm x \rangle_{\bm \mu}^+ := \mathrm{col}_i\{\left\langle x_i\right\rangle_{\mu_i}^+\}$.
\section{Modeling of \begin{neu}{the}\end{neu} Physical Network}\label{ch:woc-physical-model}
The underlying physical network is represented by a directed graph $\mathscr G_p = (\mathcal V, \mathcal E_p)$, where the set of nodes $\mathcal V = \mathcal V_{\mathcal G} \cup \mathcal V_{\mathcal I} \cup \mathcal V_{\mathcal L}$ is partitioned into generator ($\mathcal G$), inverter ($\mathcal I$), and load nodes ($\mathcal L$).
$\mathcal G$ nodes are fed by SGs of e.g. gas or hydro turbines.
$\mathcal I$ nodes are connected to power electronics interfaced RESs such as photovoltaic power stations and $\mathcal L$ nodes are connected to in\begin{neu}elastic\end{neu} consumers only.
\begin{figure}
\centering
\includegraphics[width=\columnwidth]{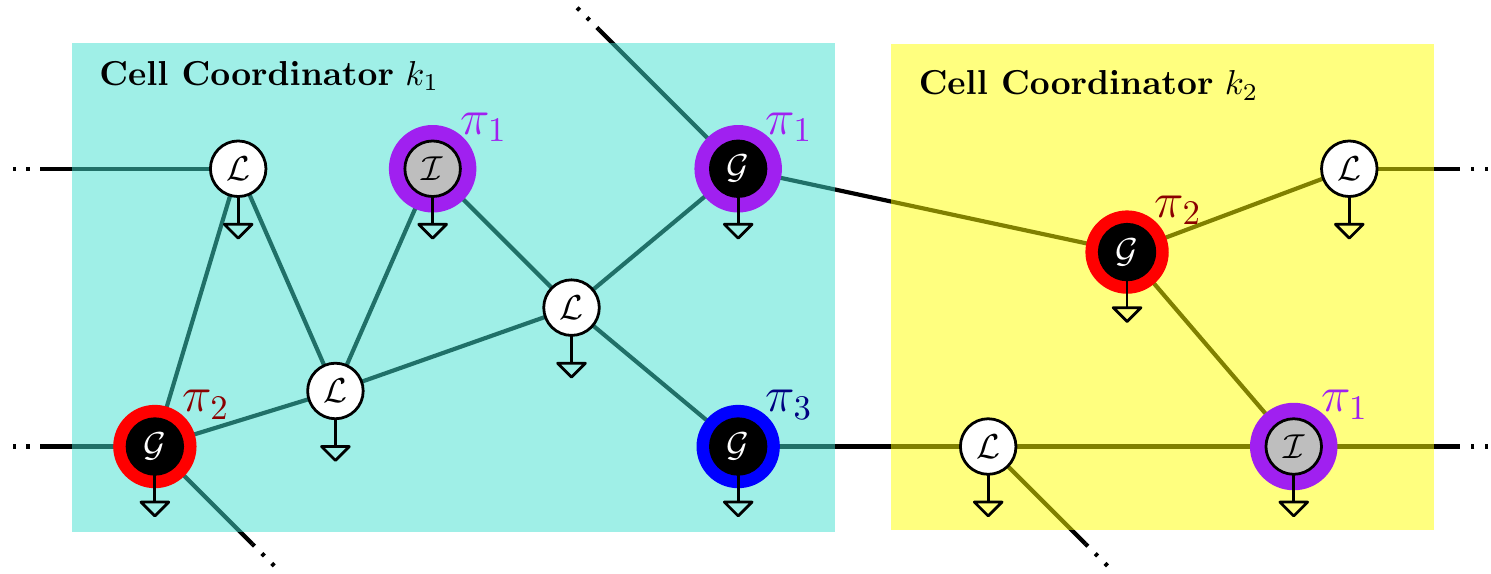}
\caption{Schematic WoC network with $\mathcal P = \{\pi_1, \pi_2, \pi_3\}$ and $\mathcal Z=\{k_1, k_2\}$. Each node of the network is of either $\mathcal G$ (black), $\mathcal I$ (grey), or $\mathcal L$ type (white).}\label{fig:WOC-scheme}
\end{figure}
$\mathcal G$ and $\mathcal I$ nodes are equipped with an active power generation input $p_{g,i}$.
In\begin{neu}elastic\end{neu} consumers may also exist at $\mathcal G$, $\mathcal I$, and $\mathcal L$ nodes, which is expressed by an active power consumption input $p_{\ell,i}$ for all $i \in \mathcal V$.
Note that 
$p_{\ell,i}$ may be negative.
\begin{neu}
In accordance with \cite{Stegink.2017,Stegink.2017c,Monshizadeh.2017,Kolsch.2020}, we will rely on the following assumptions:
\begin{assumption}\label{ass:anfang-modell}
	The network is a balanced three-phased AC system, {in which the lines can be described} by single-phase $\Pi$-equivalent circuits.
\end{assumption}
\begin{assumption}
	The network is operating around the nominal frequency.
\end{assumption}
\begin{assumption}
	Subtransient dynamics of the SGs can be neglected.
\end{assumption}
\begin{assumption}\label{ass:end-modell}
The inverter interface at $\mathcal I$ nodes is equipped with an internal 
matching
controller \cite{Jouini.2016,Monshizadeh.2017}, which has fast dynamics compared to the real-time pricing controller.
\end{assumption}
Under Assumptions \ref{ass:anfang-modell}--\ref{ass:end-modell}, the following dynamic network model can be applied:
\end{neu}
\begin{align}
\dot \vartheta_{ij} &= \omega_i- \omega_j, && i,j \in \mathcal V, \label{eq:physical-grid-begin}\\
\dot L_i &= - A_i \omega_i + p_{g,i} - p_{\ell, i} - p_i, && i \in \mathcal V_{\mathcal G} \cup \mathcal V_{\mathcal I}, \label{eq:mg-model-2}\\
\tau_{d,i} \dot U_i &= U_{f,i} - U_i - ({X_{d,i}-X_{d,i}'}){U_i^{-1}} q_i, && i \in \mathcal V_{\mathcal G}, \label{eq:mg-model-3}\\
0 &= -A_i \omega_i - p_{\ell, i} - p_i , && i \in \mathcal V_{\mathcal L}, \label{eq:mg-model-4}\\
0 &= -q_{\ell,i} - q_i, && i \in \mathcal V_{\mathcal L},\label{eq:mg-model-5}
\end{align}
where $\vartheta_{ij}$ is the voltage angle deviation between nodes $i$ and $j$, $\omega_i$ is the frequency deviation from the nominal value, $L_i$ is the angular momentum deviation, $A_i$ is a positive damping coefficient,
$\tau_{d,i}$ is the direct-axis transient open-circuit time constant,
$U_{f,i}$ and $U_i$ are the excitation voltage and internal voltage magnitudes, $X_{d,i}-X_{d,i}'$ is the d-axis synchronous reactance minus transient reactance, and $q_{\ell,i}$ is the reactive power consumption.
Eq. \eqref{eq:mg-model-2} describes the swing equation model of active power exchange between node $i$ and neighboring\footnote{The set $\mathcal N_i$ contains all adjacent nodes of $i \in \mathcal V$ without consideration of the edge direction.} nodes $j \in \mathcal N_i$.
With the internal matching controller from \cite{Jouini.2016,Monshizadeh.2017}, the dynamics of the $\mathcal I$ nodes can be described in an analogous way as change of the angular momentum $L_i$ of a virtual oscillator with a virtual damping coefficient $A_i>0$.
Eq. \eqref{eq:mg-model-3} models the transient dynamics of voltage magnitudes $U_i$ of $\mathcal G$ nodes,
and \eqref{eq:mg-model-4}--\eqref{eq:mg-model-5} describe the active and reactive power conservation at $\mathcal L$ nodes. The sending-end active and reactive power flows $p_i$ and $q_i$ in \eqref{eq:mg-model-2}--\eqref{eq:mg-model-5} are evaluated by the lossy AC power flow equations \cite{Machowski.2012}
\begin{align}
p_i &= \sum_{j \in \mathcal N_i} B_{ij} U_i U_j \sin(\vartheta_{ij}) + G_{ii} U_i^2 && \nonumber \\
&+ \sum_{j \in \mathcal N_i} G_{ij} U_i U_j \cos(\vartheta_{ij}), && i \in \mathcal V, \\
q_i &= -\sum_{j \in \mathcal N_i} B_{ij} U_i U_j \cos(\vartheta_{ij}) + B_{ii} U_i^2 && \nonumber \\
&+\sum_{j \in \mathcal N_i} G_{ij} U_i U_j \sin(\vartheta_{ij}), && i \in \mathcal V \label{eq:physical-grid-end},
\end{align}
where $B_{ij}$ [$G_{ij}$] denotes the negative of the susceptance [conductance] of line $(i,j) \in \mathcal E_p$,
$B_{ii} = -b_i - \sum_{j \in \mathcal N_i} B_{ij}$
denotes the self-susceptance at node $i \in \mathcal V$, and
 $G_{ii} = -g_i - \sum_{j \in \mathcal N_i} G_{ij}$ 
denotes the self-conductance at node $i \in \mathcal V$. 

Within the model, the active power generations $p_{g,i}$, voltage magnitudes $U_{i}$ of $\mathcal I$ nodes and excitation voltages $U_{f,i}$ of $\mathcal G$ nodes constitute the controlled input variables, whereas the active and reactive power consumptions $p_{\ell,i}$ and $q_{\ell, i}$ constitute the uncontrollable (i.e. disturbance) inputs of the network. 
\section{Optimization Problems of Network Participants}\label{ch:woc-problemformulation}
Let $\mathcal P = \{1, \ldots, |\mathcal P|\}$ denote the set of PPOs and $\mathcal Z = \{1, \ldots, |\mathcal Z|\}$ be the set of CCs.
As illustrated in Fig. \ref{fig:WOC-scheme}, each node is associated to exactly one CC and each $\mathcal G$ and $\mathcal I$ node is associated to exactly one PPO. 
We note $i \in \mathcal V_\pi^{\mathcal P}$ if node $i$ belongs to PPO $\pi \in \mathcal P$ and $i \in \mathcal V_k^{\mathcal Z}$ if node $i$ belongs to CC 
$k\in \mathcal Z$.
\subsection{Producers}
Each PPO $\pi \in \mathcal P$ seeks to maximize its overall profit
\begin{align}
\mathtt P_\pi = \sum_{i \in \mathcal V_{\pi}^{\mathcal P}}\mathtt P_{\pi,i} =  \sum_{i \in \mathcal V_\pi^{\mathcal P}} -\mathtt C_i (p_{g,i}) +  \lambda_i \cdot p_{g,i} - \omega_i \cdot p_{g,i}, \label{eq:PPO-profit}
\end{align}
where 
$\mathtt C_i(p_{g,i})$ is assumed to be a convex function representing the production costs
	 at node $i$, $\lambda_i$ is the local price for the active power generation $p_{g,i}$, and $-\omega_i \cdot p_{g,i}$ is an additive penalty for active power generation during an overfrequency period (or, respectively, an additive reward for active power generation during an underfrequency period).
\begin{neu}
	In order to encourage PPOs to produce at marginal cost, we impose the following key assumption:
	\begin{assumption}\label{ass:no-market-power}
		No PPO is large enough to exercise market power.
	\end{assumption}
\begin{remark}
Assumption \ref{ass:no-market-power} is justified if the overall number of PPOs is high \cite{Reinisch.2006} and if no PPO has a substantial market share in a specific cell \cite{Staudt.2019,Sarfati.2019}.
\end{remark} \end{neu} 
	
Denote by $\bm p_{g,\pi}$, $\bm U_{f, \pi}$, and $\bm U_{\mathcal I,\pi}$ the corresponding  
vectors of $p_{g,i}$, $U_{f,i}$, and $U_{\mathcal I,i}$ associated to PPO $\pi \in \mathcal P$. 
%
%
%
%
%
The individual profit $\mathtt P_\pi$ in \eqref{eq:PPO-profit} is \begin{neu}then\end{neu}
 driven by the local prices $\lambda_i$, which are offered by the respective CC at node $i \in \mathcal V_\pi^{\mathcal P}$.
In addition, the flexibilities $\bm p_{g,\pi}$, $\bm U_{f,\pi}$, and $\bm U_{\mathcal I,\pi}$ may be bounded from above or below\footnote{It was shown in \cite[Proposition~1]{Kolsch.2020} that upper [lower] bounds of $U_i$ for $\mathcal G$ nodes can be transformed into upper [lower] bounds of $U_{f,i}$ by a linear mapping.}.
\begin{neu}With Assumption \ref{ass:no-market-power},\end{neu} each PPO $\pi \in \mathcal P$ thus aims at solving the constrained optimization problem
\begin{subequations}\label{eq:PPO-OP}
	\begin{IEEEeqnarray}{Cll}
		\max_{\bm p_{g,\pi}, \bm U_{f, \pi}, \bm U_{\mathcal I,\pi}}& \mathtt P_{\pi} & \label{eq:PPO-OP-objective}	\\
		\quad \text{s.t.}& \underline p_{g,i} \leq p_{g,i} \leq \overline p_{g,i}, & \quad i \in \mathcal V_\pi^{\mathcal P},\label{eq:PPO-OP-constraint-begin}\\ 
		& \underline U_{f,i} \leq U_{f,i} \leq \overline U_{f,i}  & \quad i \in \mathcal V_\pi^{\mathcal P} \cap \mathcal V_{\mathcal G},\label{weitere-spannungsgleichung}\\
		& \underline U_{i} \leq U_i \leq \overline U_i, & \quad i \in \mathcal V_\pi^{\mathcal P} \cap \mathcal V_{\mathcal I}\label{eq:PPO-OP-constraint-end}.
	\end{IEEEeqnarray}
\end{subequations}
\begin{neu}
	Constraint \eqref{eq:PPO-OP-constraint-begin} represents the remaining power generation capacity of PPOs that is not contracted otherwise, \end{neu}\begin{neu}e.g. on the day-ahead market or by long-term agreements.
	 The bounds in \eqref{weitere-spannungsgleichung} and \eqref{eq:PPO-OP-constraint-end} are typically imposed by the local grid code.
\end{neu}
\subsection{Consumers}
\begin{neu}In\begin{neu}elastic\end{neu} consumers at node \end{neu}$i \in \mathcal V$ are characterized by a fixed active power consumption $p_{\ell, i}$. The cumulative cost for consumers in cell $k \in \mathcal Z$ is thus
\begin{align}
\mathtt C^\ell_k = \Lambda_k \cdot (\Phi_k + \sum_{i \in \mathcal V_k^{\mathcal Z}} p_{\ell, i}), \label{eq:consumer-cost}
\end{align}
where $\Lambda_k$ denotes the cell-specific electricity price and 
\begin{align}
\Phi_k = \sum_{i \in \mathcal V_k^{\mathcal Z}} \Bigg(G_{ii} U_i^2 + \sum_{j \in \mathcal N_i} G_{ij} U_i U_j \cos(\vartheta_{ij})\Bigg)
\end{align}
represents the resistive transmission losses within cell $k$.
\begin{neu}
\begin{remark}
Note that \begin{neu}elastic\end{neu} consumers at node $i\in \mathcal V$ which aim to maximize their profits 
\begin{align}\label{eq:flexible-consumer-problem}
\mathtt P^\ell_i = \mathtt U_i(p_{\ell,i})- \lambda_i \cdot p_{\ell,i} + \omega_i \cdot p_{\ell,i},
\end{align}
with $\mathtt U_i(p_{\ell,i})$ being a concave utility function,
can be modeled as producers with negative generation $p_{g,i}=-p_{\ell,i}\leq 0$,
since the profit maximization problem following from \eqref{eq:flexible-consumer-problem} is structurally 
identical to \eqref{eq:PPO-OP}.
Hence, without loss of generality, the term \emph{consumer} will be used synonymously with in\begin{neu}elastic\end{neu} consumer for the remainder of this paper. 
\end{remark}
\end{neu}
\subsection{Cell Coordinators}
For all $i,j \in \mathcal V_k^{\mathcal Z}$, the CC $k \in \mathcal Z$ seeks to choose the
cell-specific price $\Lambda_k=\lambda_i = \lambda_j$ in such a way that its own profit $\mathtt P_k$ (revenues from consumers minus payments to the PPOs) is maximized. Consequently, with
\begin{align}
\mathtt P_k = \Lambda_k \cdot \Bigg(  \Phi_k+ \sum_{i \in \mathcal V_k^{\mathcal Z}} p_{\ell, i}\Bigg) - \Lambda_k \cdot \sum_{\substack{i \in \mathcal V_k^{\mathcal Z} \\ i \notin \mathcal V_{\mathcal L}}} p_{g,i}, \label{CC-utility}
\end{align}
each CC $k \in \mathcal Z$ aims to solve the optimization problem
	\begin{IEEEeqnarray}{C?ll}
		\quad \max_{\Lambda_k}& \mathtt P_{k}. & \label{eq:CC-OP}
	\end{IEEEeqnarray}
Following the lines of \cite{Stegink.2017}, we can use a distributed reformulation of \eqref{eq:CC-OP} as follows:
\begin{subequations}\label{eq:CC-OP-dist}
	\begin{IEEEeqnarray}{C?ll}
		\max_{\bm \lambda_k}& \sum_{i \in \mathcal V_k^{\mathcal Z}} \Big( \lambda_i \cdot \Big(  \varphi_i+  p_{\ell, i}\Big)\Big) - \sum_{\substack{i \in \mathcal V_k^{\mathcal Z} \\ i \notin \mathcal V_{\mathcal L}}} \lambda_i p_{g,i} & \label{eq:CC-OP-dist-objective} \\
		\text{s.t.} &  (\bm D_c^k)^\top \bm \lambda_k = \nullvec{}, \label{eq:CC-OP-dist-constraint}
	\end{IEEEeqnarray}
\end{subequations}
where $\bm \lambda_k = \mathrm{col}_{i \in \mathcal V_k^{\mathcal Z}}\{\lambda_i\}$ and 
\begin{align}
\varphi_i = G_{ii}U_i^2 + \sum_{j \in \mathcal N_i} G_{ij}U_i U_j \cos(\vartheta_{ij})\label{eq-klein-phi}
\end{align}
and $\bm D_c^k$ is the incidence matrix of a connected communication graph $\mathscr G_c^k = (\mathcal V_k^{\mathcal Z}, \mathcal E_{c}^k)$. 
The reformulation \eqref{eq:CC-OP-dist} of \eqref{eq:CC-OP} is exact since $\Phi_k = \sum_{i \in \mathcal V_k^{\mathcal Z}} \varphi_i$ and since \eqref{eq:CC-OP-dist-constraint} holds if and only if for each cell $k \in \mathcal Z$, all nodal prices $\lambda_i$, $\lambda_j$ with $i,j \in \mathcal V_k^{\mathcal Z}$ are equal to a cell-specific price $\Lambda_k$.
\subsection{Primal-Dual Gradient Controller}
We now derive control laws for each network participant based on the optimization problems \begin{neu}{\eqref{eq:PPO-OP} and \eqref{eq:CC-OP-dist}}\end{neu}. For PPO $\pi \in \mathcal P$, consider the Lagrangian of \eqref{eq:PPO-OP},
\begin{align*}
\mathscr L_\pi &= - \mathtt P_\pi + \bm \mu_{p-,\pi}^\top (\underline{\bm p}_{g,\pi} - \bm p_{g,\pi}) + \bm \mu_{p+,\pi}^\top ({\bm p}_{g,\pi} - \underline{\bm p}_{g,\pi})  \nonumber \\
&+ \bm \mu_{\mathcal G-,\pi}^\top (\underline{\bm U}_{f,\pi} - \bm U_{f,\pi}) + \bm \mu_{\mathcal G+,\pi}^\top ({\bm U}_{f,\pi} - \underline{\bm U}_{f,\pi})  \nonumber \\
&+ \bm \mu_{\mathcal I-,\pi}^\top (\underline{\bm U}_{\mathcal I,\pi} - \bm U_{\mathcal I,\pi}) + \bm \mu_{\mathcal I+,\pi}^\top ({\bm U}_{\mathcal I,\pi} - \underline{\bm U}_{\mathcal I,\pi}),
\end{align*}
where $\bm \mu_{(\cdot),\pi}$ denotes the vectors of Lagrange multipliers for the inequality constraints \eqref{eq:PPO-OP-constraint-begin}--\eqref{eq:PPO-OP-constraint-end}.
Since \eqref{eq:PPO-OP} is convex and Slater's Constraint Qualification \cite{Boyd.2015} is fulfilled, the KKT conditions specifying a saddle point of $\mathscr L_\pi$ can be applied to derive a necessary and sufficient condition for an optimizer of \eqref{eq:PPO-OP}:
\begin{align}
\nullvec{} & = \nabla \mathtt C_\pi(\bm p^\star_{g,\pi}) - \bm \lambda_{\pi} + \bm \omega_\pi {-\bm{\mu}^\star_{g-,\pi}+\bm{\mu}^\star_{g+,\pi}}, \\
\nullvec{} &= - \bm \mu_{\mathcal G-,\pi}^\star + \bm \mu_{\mathcal G+,\pi}^\star,\\
\nullvec{} &=  \bm \mu_{\mathcal I-,\pi}^\star  - \bm \mu_{\mathcal I+,\pi}^\star, \\
\nullvec{} & = (\bm \mu_{g-,\pi}^\star)^\top (\underline{\bm p}_{g,\pi} - {{\bm p}^\star_{g,\pi}}), \\
\nullvec{} & =(\bm \mu_{g+,\pi}^\star)^\top( {\bm p}_{g,\pi}^\star - \overline{\bm p}_{g,\pi} ), \\
\nullvec{} &= (\bm \mu_{\mathcal G-,\pi}^\star)^\top \left(\underline{\bm U}_{f,\pi} - \bm U_{f,\pi}^\star \right), \\
\nullvec{} & = (\bm \mu_{\mathcal G+,\pi}^\star)^\top \left(  \bm U^\star_{f,\pi}   - \overline{\bm U}_{f,\pi} \right), \\
\nullvec{} &=(\bm \mu_{\mathcal I-,\pi}^\star)^\top \left( \underline{\bm U}_{\mathcal I,\pi} - \bm U^\star_{\mathcal I,\pi} \right),\\
\nullvec{} & = (\bm \mu_{\mathcal I+,\pi}^\star)^\top \left(\bm U^\star_{\mathcal I,\pi} - \overline{\bm U}_{\mathcal I,\pi} \right),  \\
\nullvec{} & \leq \bm \mu^\star_{\mathcal G-,\pi}, \bm \mu^\star_{\mathcal G+,\pi}, \bm \mu^\star_{\mathcal I-,\pi}, \bm \mu^\star_{\mathcal I+,\pi}, \bm \mu^\star_{g-,\pi}, \bm \mu^\star_{g+,\pi},
\end{align}
where $\nabla \mathtt C_\pi(\bm p_{g,\pi}) = \mathrm{diag}_{i \in \mathcal V_\pi^{\mathcal P}}\{\nabla \mathtt C_i(p_{g,i})\}$.
Application of the primal-dual gradient method \cite{Arrow.1958} leads to the following controller equations:
\begin{align}
\bm \tau_{g,\pi} \dot{\bm p}_{g,\pi} &= 
- \nabla \mathtt C_\pi(\bm p_{g,\pi}) + {\bm \lambda}_\pi - \bm \omega_\pi
+\bm{\mu}_{g-,\pi}-\bm{\mu}_{g+,\pi},\label{eq:PPO-controller-begin}\\
\bm \tau_{\bm U_{\mathcal G,\pi}} \dot{\bm U}_{f,\pi}&=  \bm \mu_{\mathcal G-,\pi} - \bm \mu_{\mathcal G+,\pi},\\
\bm \tau_{\bm U_{\mathcal I,\pi}} \dot{\bm U}_{\mathcal I,\pi}&=  
- \bm \mu_{\mathcal I-,\pi}  + \bm \mu_{\mathcal I+,\pi},\\
{\bm \tau_{\bm{\mu}_{g-,\pi}}\dot{\bm{\mu}}_{g-,\pi}} & = {\langle \underline{\bm p}_{g,\pi} - \bm p_{g,\pi} \rangle^+_{\bm{\mu}_{g-,\pi}}},\\
{\bm \tau_{\bm{\mu}_{g+,\pi}}\dot{\bm{\mu}}_{g+,\pi}} &= {\langle  \bm p_{g,\pi} - \overline{\bm p}_{g,\pi}\rangle^+_{\bm{\mu}_{g+,\pi}}}, \\
\bm \tau_{\bm{\mu}_{\mathcal G-,\pi}}\dot{\bm{\mu}}_{\mathcal G-,\pi} &= \langle \underline{\bm{U}}_{f,\pi} - \bm U_{f,\pi} \rangle^+_{\bm{\mu}_{\mathcal G-,\pi}}, \\
\bm \tau_{\bm{\mu}_{\mathcal G+,\pi}}\dot{\bm{\mu}}_{\mathcal G+,\pi} &= \langle  \bm U_{f,\pi} - \overline{\bm{U}}_{f,\pi}\rangle^+_{\bm{\mu}_{\mathcal G+,\pi}}, \\
\bm \tau_{\bm{\mu}_{\mathcal I-,\pi}}\dot{\bm{\mu}}_{\mathcal I-,\pi} &= \langle \underline{\bm U}_{\mathcal I,\pi} - \bm U_{\mathcal I,\pi} \rangle^+_{\bm{\mu}_{\mathcal I-,\pi}},\\
\bm \tau_{\bm{\mu}_{\mathcal I+,\pi}}\dot{\bm{\mu}}_{\mathcal I+,\pi} &= \langle  \bm U_{\mathcal I,\pi} - \overline{\bm{U}}_{\mathcal I,\pi}\rangle^+_{\bm{\mu}_{\mathcal I+,\pi}}.\label{eq:PPO-controller-end}
\end{align}
with diagonal matrices $\bm \tau_{(\cdot)}\succ 0$.

\begin{neu}{For CC $k \in \mathcal Z$,}\end{neu}
the Lagrangian of \eqref{eq:CC-OP-dist} equals
\begin{align*}
\mathscr L_k = - \sum_{i \in \mathcal V_k^{\mathcal Z}} ( \lambda_i \cdot (  \varphi_i+  p_{\ell, i})) + \sum_{\substack{i \in \mathcal V_k^{\mathcal Z} \\ i \notin \mathcal V_{\mathcal L}}} \lambda_i p_{g,i} + \bm \nu_k^\top (\bm D_c^k)^\top \bm \lambda_k,
\end{align*}
where $\bm \nu_k$ is the Lagrange multiplier associated to the linear equality constraint \eqref{eq:CC-OP-dist-constraint}.
Since \eqref{eq:CC-OP-dist} is convex and Slater's Constraint Qualification is fulfilled, we can again specify the KKT conditions
\begin{align}
\nullvec{} &= \bm p_{g,k} - \bm \varphi_k - \bm p_{\ell,k} + \bm D_c^k \bm \nu_k^\star, \label{eq:CC-constr1}\\
\nullvec{} & = (\bm D_c^k)^\top \bm \lambda_k^\star.  \label{eq:CC-constr2}
\end{align}
Application of the primal-dual gradient method yields
\begin{align}
\bm \tau_{\bm \lambda,k} \dot{\bm \lambda}_k&= -\bm p_{g,k} + \bm \varphi_k + \bm p_{\ell,k} - \bm D_c^k \bm \nu_k, \label{eq:CC-controller-begin}\\
\bm \tau_{\bm \nu,k} \dot{\bm \nu}_k& = (\bm D_c^k)^\top \bm \lambda_k.\label{eq:CC-controller-end}
\end{align}
The simultaneous execution of the optimization schemes of PPOs $\pi \in \mathcal P$ and CCs $k \in \mathcal Z$ leads to a superposition of the derived control equations \eqref{eq:PPO-controller-begin}--\eqref{eq:PPO-controller-end}, \eqref{eq:CC-controller-begin}--\eqref{eq:CC-controller-end}.
Consequently, the resulting closed-loop system consists of the physical system \eqref{eq:physical-grid-begin}--\eqref{eq:physical-grid-end} along with the controller equations
\begin{align}
\bm \tau_{g} \dot{\bm p}_{g} &= 
- \nabla \mathtt C(\bm p_{g}) + {\bm \lambda} - \bm \omega
+\bm{\mu}_{g-}-\bm{\mu}_{g+},\label{eq:WOC-closed-loop-begin}\\
\bm \tau_{\bm U_{\mathcal G}} \dot{\bm U}_{f}&=  \bm \mu_{\mathcal G-} - \bm \mu_{\mathcal G+},\label{eq:WOC-closed-loop-zweites}\\
\bm \tau_{\bm U_{\mathcal I}} \dot{\bm U}_{\mathcal I}&=  
- \bm \mu_{\mathcal I-}  + \bm \mu_{\mathcal I+},\\
{\bm \tau_{\bm{\mu}_{g-}}\dot{\bm{\mu}}_{g-}} & = {\langle \underline{\bm p}_{g} - \bm p_{g} \rangle^+_{\bm{\mu}_{g-}}},\\
{\bm \tau_{\bm{\mu}_{g+}}\dot{\bm{\mu}}_{g+}} &= {\langle  \bm p_{g} - \overline{\bm p}_{g}\rangle^+_{\bm{\mu}_{g+}}}, \\
\bm \tau_{\bm{\mu}_{\mathcal G-}}\dot{\bm{\mu}}_{\mathcal G-} &= \langle \underline{\bm{U}}_{f} - \bm U_{f} \rangle^+_{\bm{\mu}_{\mathcal G-}}, \\
\bm \tau_{\bm{\mu}_{\mathcal G+}}\dot{\bm{\mu}}_{\mathcal G+} &= \langle  \bm U_{f} - \overline{\bm{U}}_{f}\rangle^+_{\bm{\mu}_{\mathcal G+}}, \\
\bm \tau_{\bm{\mu}_{\mathcal I-}}\dot{\bm{\mu}}_{\mathcal I-} &= \langle \underline{\bm U}_{\mathcal I} - \bm U_{\mathcal I} \rangle^+_{\bm{\mu}_{\mathcal I-}},\\
\bm \tau_{\bm{\mu}_{\mathcal I+}}\dot{\bm{\mu}}_{\mathcal I+} &= \langle  \bm U_{\mathcal I} - \overline{\bm{U}}_{\mathcal I}\rangle^+_{\bm{\mu}_{\mathcal I+}},\label{eq:WOC-closed-loop-PPO-end}\\
\bm \tau_{\bm \lambda} \dot{\bm \lambda}&= -\bm p_{g} + \bm \varphi + \bm p_{\ell} - \bm D_c \bm \nu, \label{eq:WOC-closed-loop-begin-CC}\\
\bm \tau_{\bm \nu} \dot{\bm \nu}& = \bm D_c^\top \bm \lambda,\label{eq:WOC-closed-loop-end}
\end{align}
where the parameters and variables in \eqref{eq:WOC-closed-loop-begin}--\eqref{eq:WOC-closed-loop-end} are stacked vectors or diagonal matrices of appropriate sizes for all $\pi \in \mathcal P$ and $k \in \mathcal Z$, e.g. $\bm \tau_g = \mathrm{diag}_{\pi \in \mathcal P}\{\bm \tau_{g,\pi}\}$ and $\bm D_c = \mathrm{col}_{k \in \mathcal Z}\{\bm D_c^k\}$, and $\mathtt C(\bm p_g)$ equals the sum of all cost functions $\mathtt C_\pi(\bm p_{g,\pi})$ in \eqref{eq:PPO-profit}.
\begin{remark}\label{rem:Dcblockdiagonal}
	Note that $\bm D_c$ has a block diagonal structure, thus the corresponding communication graph $\mathscr G_c$ is disconnected and contains $|\mathcal Z|$ connected subgraphs.
\end{remark}
Before discussing the interactions of cell-specific prices, we point out a key result concerning the zero-frequency deviation of the closed-loop system \eqref{eq:physical-grid-begin}--\eqref{eq:physical-grid-end}, \eqref{eq:WOC-closed-loop-begin}--\eqref{eq:WOC-closed-loop-end}:
\begin{proposition}[Zero Frequency Deviation]\label{prop:frequenz-null-WOC}
	At each equilibrium of
 \eqref{eq:physical-grid-begin}--\eqref{eq:physical-grid-end}, \eqref{eq:WOC-closed-loop-begin}--\eqref{eq:WOC-closed-loop-end}, it holds that $\omega_i=0$ for all $i\in \mathcal V$, i.e. 
 each node is operating at the nominal frequency.
\end{proposition}
\begin{IEEEproof}
	Eq. \eqref{eq:physical-grid-begin} can be written in vector-matrix notation as $\dot{\bm \vartheta} = \bm D_p^\top \bm \omega$, where $\bm D_p$ denotes the incidence matrix of the physical network $\mathscr G_p$. Since $\mathscr G_p$ is a connected graph, each equilibrium $\bm \omega^\star$ with $\nullvec{} = \bm D_p^\top \bm \omega^\star$ fulfills $\bm \omega^\star = \omega^\star \cdot \mathbb 1$, i.e. the nodal frequencies are synchronized at steady state. 
	Since $\mathbb 1^\top \bm \varphi = \Phi$ and $\mathbb 1^\top \bm D_c = \nullvec{}$, left-multiplying \eqref{eq:WOC-closed-loop-begin-CC} by $\mathbb 1^\top$ and inserting the equilibrium values implies
	\begin{align}
	0 = - \sum_{i \in \mathcal V_{\mathcal G} \cup \mathcal V_{\mathcal I}} p_{g,i}^\star + \sum_{i \in \mathcal V} p_{\ell, i} + \Phi.\label{eq:proof-freq-0-1}
	\end{align}
	Due to the fact that $\Phi = \sum_{i \in \mathcal V} p_i$, a summation of all equations in \eqref{eq:mg-model-2} and \eqref{eq:mg-model-4} and insertion of the equilibrium values leads to
	\begin{align}
	0 = -\sum_{i \in \mathcal V_{\mathcal G} \cup \mathcal V_{\mathcal I}} A_i  \omega_i^\star + \sum_{i \in \mathcal V_{\mathcal G} \cup \mathcal V_{\mathcal I}} p_{g,i}^\star -  \sum_{i \in \mathcal V} p_{\ell, i} - \Phi. \label{eq:proof-freq-0-2}
	\end{align}
	Comparison between \eqref{eq:proof-freq-0-1} and \eqref{eq:proof-freq-0-2} yields
	\begin{align}
	0 = -\sum_{i \in \mathcal V_{\mathcal G} \cup \mathcal V_{\mathcal I}} A_i \omega_i^\star = - \omega^\star \cdot \sum_{i \in \mathcal V_{\mathcal G} \cup \mathcal V_{\mathcal I}} A_i.
	\end{align}
	Finally, since $A_i>0$ holds by definition, it follows that $\omega^\star$ is zero.
\end{IEEEproof}
\section{Interaction Between Energy Cells}\label{ch:woc-multiple-cells}
Equation \eqref{eq:WOC-closed-loop-end} implies that at each equilibrium of 
\eqref{eq:WOC-closed-loop-begin}--\eqref{eq:WOC-closed-loop-end},
nodal prices $\lambda_i, \lambda_j \in \mathcal V_k^{\mathcal Z}$ are equal to a common zonal price $\Lambda_k$ for each cell $k \in \mathcal Z$.
Incentives for increased or decreased power generation can thus be imposed by means of differences in zonal prices.
However, Remark \ref{rem:Dcblockdiagonal} implies that there is no direct relationship between each individual $\Lambda_k$.
\begin{neu}{In the next subsection, we introduce specific couplings between CCs to enable such relationships.}\end{neu}

\subsection{Coupling of Zonal Prices $\Lambda_k$}
To provide interdependencies between the zonal prices, additional constraints of the form
\begin{align}
\lambda_i \stackrel != \eta_{ij} \cdot \lambda_j, \label{eq:kopplungs-constraints}
\end{align}
may be imposed for pairs of nodes $i,j \in \mathcal V$, which are located in \emph{different} cells, where $\eta_{ij}>0$ is an appropriate multiplier.
If \eqref{eq:CC-OP-dist-constraint} is extended by these additional constraints \eqref{eq:kopplungs-constraints}, both \eqref{eq:CC-OP-dist-constraint} and \eqref{eq:kopplungs-constraints} can be combined to the extended constraint 
\begin{align}
\nullvec{} = (\bm D_c^+)^\top \bm \lambda.\label{eq:extended-constraint-D+}
\end{align} 
Hence, the resulting controller equations 
\eqref{eq:WOC-closed-loop-begin-CC}--\eqref{eq:WOC-closed-loop-end} become
\begin{align}
\bm \tau_{\bm \lambda} \dot{\bm \lambda}&= -\bm p_{g} + \bm \varphi + \bm p_{\ell} - \bm D_c^+ \bm \nu, \label{eq:WOC-closed-loop-Dplus-begin}\\
\bm \tau_{\bm \nu} \dot{\bm \nu}& = (\bm D_c^+)^\top \bm \lambda.\label{eq:WOC-closed-loop-Dplus-end}
\end{align}
With this notation, $\bm D_c^+$ can be interpreted as the incidence matrix of an extended, weighted communication graph $\mathscr G_c^+=(\mathcal V, (\mathcal E_c,\mathcal E_c^{b}))$, where $\mathcal E_c^b$ represents communication across cell boundaries.
The weights of the edges of $\mathscr G_c^+$ are equal to $\eta_{ij}$, if the two adjacent nodes belong to different cells, and equal to $1$ if $i$ and $j$ are located in the same cell. 
\begin{assumption}\label{ass:eta-feasible}
The multipliers $\eta_{ij}$ are chosen in a \emph{feasible} sense such that 
there exists at least one $\bm \lambda > \nullvec{n}$ fulfilling \eqref{eq:extended-constraint-D+}.
In this case, each cell $k \in \mathcal Z$ can be characterized with a specific participation factor $\kappa_k > 0$ such that each $\eta_{ij}$ is calculated by 
$
\eta_{ij} = {\kappa_{k_1}}/{\kappa_{k_2}}
$
if $i$ is located in cell $k_1$ and $j$ is located in cell $k_2$.
\end{assumption}
\begin{proposition}[Connectivity of Zonal Prices]
	If $\mathscr G_c^+$ is weakly connected, then
	\begin{align}
	\frac{\lambda_1}{\kappa_{k_1}} = \frac{\lambda_2}{\kappa_{k_2}} = \cdots = \frac{\lambda_n}{\kappa_{k_n}} =: \Lambda^\circ, \label{eq:es-gibt-einen-fair-price}
	\end{align}
where $\kappa_{k_i}$ denotes the corresponding participation factor belonging to the cell $k$ where node $i \in \mathcal V$ is located.
\end{proposition}
\begin{IEEEproof}
Define the auxiliary matrices
$\bm K_1 = \mathrm{diag}_{ij \in (\mathcal E_c, \mathcal E_c^b)}\{\eta_{ij}\} \succ 0$ and
$\bm K_2 = \mathrm{diag}_{i \in \mathcal V}\{\kappa_{k_i}\}  \succ 0$.
Then, $(\bm D_c^+)^\top\bm \lambda = \nullvec{}$ is equivalent to
\begin{align}
\underbrace{\bm K_1^{-1} \nullvec{}}_{\nullvec{}} = \bm K_1^{-1} (\bm D_c^+)^\top \bm \lambda = \underbrace{\bm K_1^{-1} (\bm D_c^+)^\top \bm K_2}_{=:(\bm D_c^\circ)^\top} \underbrace{\bm K_2^{-1} \bm \lambda}_{=:\bm \lambda^\circ}. \label{eq:D+projected}
\end{align}
Inserting the definition of $\bm K_1$ into \eqref{eq:D+projected} with $\eta_{ij} = {\kappa_{k_1}}/{\kappa_{k_2}}$ reveals that $\bm D_c^\circ$ is the incidence matrix of a new communication graph $\mathscr G_c^\circ$, which is equivalent to $\mathscr G_c^+$ with all edge weights reset to $1$. Thus each solution $\bm \lambda^\circ$ of the resulting equation $\nullvec{} = (\bm D_c^\circ)^\top \bm \lambda^\circ$ is of the form $\bm \lambda^\circ = \mathds 1 \cdot \mathrm{const}$, i.e. each component of $\bm \lambda^\circ$ has the same value.
Finally, since $\bm \lambda^\circ = \bm K_2^{-1}\bm \lambda =  \mathrm{col}_{i\in \mathcal V}\{\lambda_i / \kappa_{k_i}\}$, this leads to \eqref{eq:es-gibt-einen-fair-price}.
\end{IEEEproof}	
\begin{remark}
	If $\mathscr G_c^+$ is not weakly connected, then \eqref{eq:es-gibt-einen-fair-price} holds separately for all nodes in each weakly connected component in  $\mathscr G_c^+$.
\end{remark}
\begin{remark}
The uniform price $\Lambda^\circ$ resulting when $\bm \kappa=\mathds 1$ is called the \emph{market-clearing price} of the network, since in this case the total revenue of all PPOs is equal to the cumulative costs of all CCs and consumers.
Accordingly, $\kappa_k$ describes the multiplicity of the cell-specific price $\Lambda_k$ compared to $\Lambda^\circ$.
\end{remark}
In the following, we will discuss to what extent $\bm \kappa$ can be applied to modify the impact of specific cells on the overall network.
For this purpose, we first investigate the equivalence of the multiple optimization problems of PPOs and CCs to a modified, centralized optimization problem.
\subsection{Comparison with Centralized Optimization}
\begin{proposition}[Equivalence to Centralized Optimization]\label{prop:equivalence-to-centralized-optimization}
	Define the centralized optimization problem
	\begin{subequations}\label{eq:central-OP}
		\begin{IEEEeqnarray}{Cll}
			\quad \max_{\bm p_{g}, \bm U_{f}, \bm U_{\mathcal I}}& \mathtt P^\kappa \label{eq:central-OP-objective}	\\
			\quad \text{s.t.} & \bm \Phi = \sum_{i \in \mathcal V_{\mathcal G} \cup \mathcal V_{\mathcal I}} p_{g,i} - \sum_{i \in \mathcal V} p_{\ell,i}, \label{eq:central-OP-constraint-begin}\\
			& \underline{\bm{p}}_{g} \leq \bm p_{g} \leq \overline{\bm{p}}_{g},\label{eq:central-OP-inequality-constraint-begin}\\ 
			& \underline{\bm{U}}_{f} \leq \bm U_{f} \leq \overline{\bm{U}}_{f}, \\
			& \underline{\bm{U}}_{\mathcal I} \leq \bm U_{\mathcal I} \leq \overline{\bm{U}}_{\mathcal I}, \label{eq:central-OP-constraint-end}
		\end{IEEEeqnarray}
	\end{subequations}
where
\begin{align}
\mathtt P^\kappa = -\sum_{k \in \mathcal Z} \sum_{i \in \mathcal V_k^{\mathcal Z}} \frac{1}{\kappa_k} \cdot \mathtt C_i(p_{g,i}) - \sum_{i\in \mathcal V_{\mathcal G} \cup \mathcal V_{\mathcal I}} \omega_i \cdot p_{g,i} .\label{eq:Ckappa-definition}
\end{align}
Then each optimizer $(\bm p_g^\star, \bm U_f^\star, \bm U_{\mathcal I}^\star)$ of \eqref{eq:central-OP} is an equilibrium of \eqref{eq:WOC-closed-loop-begin}--\eqref{eq:WOC-closed-loop-PPO-end}, \eqref{eq:WOC-closed-loop-Dplus-begin}--\eqref{eq:WOC-closed-loop-Dplus-end} and vice versa.
If $\mathtt C_i(p_{g,i})$ in \eqref{eq:Ckappa-definition} are strictly convex, then $(\bm p_g^\star, \bm U_f^\star, \bm U_{\mathcal I}^\star)$ is unique.
\end{proposition}
\begin{IEEEproof}
Constraint \eqref{eq:central-OP-constraint-begin} is equivalent to (cf. \cite[p.~2615]{Stegink.2017})
	\begin{align}
	\widehat{\bm D}_c \widehat{\bm \nu} = \bm p_g - \bm p_{\ell} - \bm \varphi, \label{eq:frequency-reformulation}
	\end{align}
where $\bm \varphi = \mathrm{col}_i\{\varphi_i\}$ (see \eqref{eq-klein-phi}) and $\widehat{\bm D}_c$ is the incidence matrix of a connected communication graph. If we choose $\widehat{\bm D}_c = \bm D_c^\circ$
and define
$\mathtt C^\kappa(\bm p_g) = \sum_{k \in \mathcal Z} \sum_{i \in \mathcal V_k^{\mathcal Z}} \frac{1}{\kappa_k} \cdot \mathtt C_i(p_{g,i})$,
then the Lagrangian of \eqref{eq:central-OP} becomes
\begin{align}
\mathscr L^\kappa &= \mathtt C^\kappa(\bm p_g)  +  \sum_{i\in \mathcal V_{\mathcal G} \cup \mathcal V_{\mathcal I}} \omega_i \cdot p_{g,i} \nonumber \\
&+ \widehat{\bm\lambda}^\top ({\bm D}_c^\circ \widehat{\bm \nu} - \bm p_g + \bm p_\ell + \bm \varphi) \nonumber \\
&+ \widehat{\bm \mu}_{g-}^\top (\underline{\bm p}_{g} - \bm p_{g}) + \widehat{\bm \mu}_{g+}^\top ({\bm p}_{g} - \underline{\bm p}_{g})  \nonumber \\
&+ \widehat{\bm \mu}_{\mathcal G-}^\top (\underline{\bm U}_{f} - \bm U_{f}) + \widehat{\bm \mu}_{\mathcal G+}^\top ({\bm U}_{f} - \underline{\bm U}_{f})  \nonumber \\
&+ \widehat{\bm \mu}_{\mathcal I-}^\top (\underline{\bm U}_{\mathcal I} - \bm U_{\mathcal I}) + \widehat{\bm \mu}_{\mathcal I+}^\top ({\bm U}_{\mathcal I} - \underline{\bm U}_{\mathcal I}),
\end{align}
where $\widehat{\bm \lambda}$ denotes the Lagrange multiplier for equality constraint \eqref{eq:frequency-reformulation} and
$\widehat{\bm \mu}_{(\cdot)}$ are the Lagrange multipliers for the inequality constraints \eqref{eq:central-OP-inequality-constraint-begin}--\eqref{eq:central-OP-constraint-end}.
Since \eqref{eq:central-OP} is convex and Slater's Constraint Qualification is fulfilled, the primal-dual optimizer of \eqref{eq:central-OP} is given by
\begin{align}
\nullvec{} &= 
- \nabla \mathtt C^\kappa(\bm p_{g}^\sharp) - \bm \omega_i + \widehat{\bm \lambda}^\sharp 
+\widehat{\bm{\mu}}_{g-}^\sharp-\widehat{\bm{\mu}}_{g+}^\sharp,\label{eq:central-closed-loop-begin}\\
\nullvec{}&=  \widehat{\bm \mu}_{\mathcal G-}^\sharp - \widehat{\bm \mu}_{\mathcal G+}^\sharp,\\
\nullvec{}&=  
- \widehat{\bm \mu}_{\mathcal I-}^\sharp  + \widehat{\bm \mu}_{\mathcal I+}^\sharp,\\
\nullvec{} & = (\widehat{\bm{\mu}}_{g-}^\sharp)^\top( \underline{\bm p}_{g} - \bm p_{g}^\sharp ),\\
\nullvec{} &= (\widehat{\bm{\mu}}_{g+}^\sharp)^\top(  \bm p_{g}^\sharp - \overline{\bm p}_{g}), \\
\nullvec{} &= (\widehat{\bm{\mu}}_{\mathcal G-}^\sharp)^\top(\underline{\bm{U}}_{f} - \bm U_{f}^\sharp), \\
\nullvec{} &= (\widehat{\bm{\mu}}_{\mathcal G+}^\sharp)^\top(  \bm U_{f}^\sharp - \overline{\bm{U}}_{f}), \\
\nullvec{} &= (\widehat{\bm{\mu}}_{\mathcal I-}^\sharp)^\top( \underline{\bm U}_{\mathcal I} - \bm U_{\mathcal I}^\sharp),\\
\nullvec{} &= (\widehat{\bm{\mu}}_{\mathcal I+}^\sharp)^\top( \bm U_{\mathcal I}^\sharp - \overline{\bm{U}}_{\mathcal I}),\\
\nullvec{}&= -\bm p_{g}^\sharp + \bm \varphi + \bm p_{\ell} - {\bm D}_c^\circ \widehat{\bm \nu}^\sharp, \\
\nullvec{}& = ({\bm D}_c^\circ)^\top \widehat{\bm \lambda}^\sharp, \label{eq:central-Dlambda}\\
\nullvec{}& \leq \bm \mu^\sharp_{\mathcal G-}, \bm \mu^\sharp_{\mathcal G+}, \bm \mu^\sharp_{\mathcal I-}, \bm \mu^\sharp_{\mathcal I+}, \bm \mu^\sharp_{g-}, \bm \mu^\sharp_{g+}.\label{eq:central-closed-loop-end}
\end{align}
Inserting the definition \eqref{eq:D+projected} in \eqref{eq:central-Dlambda} and comparing with the right-hand side of \eqref{eq:WOC-closed-loop-end} yields
$
\widehat{\bm \lambda}^\sharp = \bm K_2^{-1} \bm \lambda^\star
$.
Hence with $ \nabla \mathtt C^\kappa(\bm p_{g}^\sharp) = \bm K_2^{-1} \nabla \mathtt C(\bm p_g)$ in \eqref{eq:central-closed-loop-begin} and by comparing \eqref{eq:central-closed-loop-begin}--\eqref{eq:central-closed-loop-end} with \eqref{eq:WOC-closed-loop-begin}--\eqref{eq:WOC-closed-loop-end} we get the equivalences
$\bm p_g^\sharp = \bm p_g^\star$,
$\bm U_f^\sharp = \bm U_f^\star$,
$\bm U_{\mathcal I}^\sharp = \bm U_{\mathcal I}^\star$,
$\widehat{\bm \mu}_{g-}^\sharp = \bm K_2^{-1} \bm \mu_{g-}^\star$,
$\widehat{\bm \mu}_{g+}^\sharp = \bm K_2^{-1} \bm \mu_{g+}^\star$,
$\widehat{\bm \mu}_{\mathcal G-}^\sharp =  \bm \mu_{\mathcal G-}^\star$,
$\widehat{\bm \mu}_{\mathcal G+}^\sharp = \bm \mu_{\mathcal G+}^\star$,  
$\widehat{\bm \mu}_{\mathcal I-}^\sharp = \bm \mu_{\mathcal I-}^\star$, 
$\widehat{\bm \mu}_{\mathcal I+}^\sharp =  \bm \mu_{\mathcal I+}^\star$.
Moreover, the dual optimizer $\bm \nu^\star$ in \eqref{eq:WOC-closed-loop-begin-CC} and the primal optimizer $\widehat{\bm \nu}^\sharp$ of \eqref{eq:frequency-reformulation} are connected via the relationship
\begin{align}
-\bm D_c^\circ \bm \nu^\sharp + \bm D_c^+ \bm \nu^\sharp = \nullvec{}. \label{eq:zwei-nus-zusammenhaenge}
\end{align}
From \eqref{eq:D+projected} it follows that $\bm D_c^\circ = \bm K_2 \bm D_c^+ \bm K_1^{-1}$ with $\bm K_1, \bm K_2 \succ 0$, thus the images of $\bm D_c^+$ and $\bm D_c^\circ$ are identical. Accordingly, for each $\bm \lambda^\star$ there exists a $\widehat{\bm \lambda}^\sharp$ fulfilling \eqref{eq:zwei-nus-zusammenhaenge}, and vice versa.

In summary, for each primal-dual optimizer of \eqref{eq:central-OP} there exists exactly one corresponding equilibrium of \eqref{eq:WOC-closed-loop-begin}--\eqref{eq:WOC-closed-loop-PPO-end}, \eqref{eq:WOC-closed-loop-Dplus-begin}--\eqref{eq:WOC-closed-loop-Dplus-end}. In particular, $(\bm p_g^\sharp, \bm U_f^\sharp, \bm U_{\mathcal I}^\sharp)=(\bm p_g^\star, \bm U_f^\star, \bm U_{\mathcal I}^\star)$.
Since \eqref{eq:central-OP} is a convex optimization problem, convergence of the trajectory $(\bm p_g(t), \bm U_f(t), \bm U_{\mathcal I}(t))$ 
to $(\bm p_g^\star, \bm U_f^\star, \bm U_{\mathcal I}^\star)$ is guaranteed (cf. \cite[Theorem~2]{Antipin.1994}). 

	If the cost functions  $\mathtt C_i(p_{g,i})$ are strictly convex, then $\mathtt P^\kappa$ is strictly concave. Hence, the equilibrum $(\bm p_g^\sharp, \bm U_f^\sharp, \bm U_{\mathcal I}^\sharp)$ from centralized optimization and therewith the equilibrium $(\bm p_g^\star, \bm U_f^\star, \bm U_{\mathcal I}^\star)$ from distributed optimization are unique.
\end{IEEEproof}
\subsection{Analysis of Pareto Efficiency}
As stated in Proposition \ref{prop:equivalence-to-centralized-optimization}, the interaction of distributed PPOs and CCs, subject to cell-specific pricing by $\kappa_k$, leads to the same equilibrium $(\bm p_g, \bm U_f, \bm U_{\mathcal I})$ as if a centralized authority with full knowledge of the whole network would solve the optimization problem \eqref{eq:central-OP} with modified cost $\mathtt C^\kappa(\bm p_g)$, where the cumulative cost functions of each cell $k \in \mathcal Z$ are divided by $\kappa_k$.  
This equivalence reveals some further key properties of the distributed WoC controller that allows the evaluation of the Pareto efficiency of the equilibrium from a multi-objective perspective.
\begin{remark}
	If $\bm \kappa = \mathbb 1$, then the objective function of \eqref{eq:central-OP} is
	\begin{align}
	\mathtt P^\kappa = \sum_{\pi \in \mathcal P} \mathtt P^{\pi} - \sum_{k \in \mathcal Z} \mathtt C_k^\ell + \sum_{k \in \mathcal Z} \mathtt P_k = -\sum_{i \in \mathcal V} \mathtt C_i(p_i),
	\end{align}
	and thus equal to the sum of all payoffs of PPOs \eqref{eq:PPO-profit}, consumers \eqref{eq:consumer-cost}, and CCs \eqref{CC-utility}. Accordingly, the distributed controller \eqref{eq:WOC-closed-loop-begin}--\eqref{eq:WOC-closed-loop-PPO-end}, \eqref{eq:WOC-closed-loop-Dplus-begin}--\eqref{eq:WOC-closed-loop-Dplus-end} leads to a constrained minimization of the overall costs of power production.
\end{remark}
As mentioned earlier, the participation factor $\kappa_k$  can be used to increase (large $\kappa_k$) or decrease (small $\kappa_k$) the proportion of active power generation in cell $k \in \mathcal Z$ as compared to the other cells. 
It directly follows from Proposition \ref{prop:equivalence-to-centralized-optimization} that the result of the optimization procedure \eqref{eq:central-OP} is invariant with respect to a scalar multiplication of $\bm \kappa$, thus the level of excess or shortage of generation
in cell $k$ does not depend on the absolute value of $\kappa_k$, but on the relative value compared to the other components in $\bm \kappa$.

Besides the fact that all $\bm \kappa > \nullvec{}$ are feasible (cf. Assumption \ref{ass:eta-feasible}), it has to be stressed that each $\bm \kappa$ achieves an \emph{efficient} allocation in the sense that for a given $\bm \kappa'>0$, there is no possibility to find a ``better'' $\bm \kappa'' \neq \bm \kappa'$, $\bm \kappa'' >0$ such that at least one cell is at a lower cost and no cell is at a higher cost.
This is formalized in the next proposition.
\begin{proposition}[Pareto efficiency of $\bm \kappa$]
	Define by
	\begin{align}
	\mathtt C_k^\mathrm{total} = -\Bigg(\sum_{i \in \mathcal V_k^{\mathcal Z}} \sum_{\pi \in \mathcal P} \mathtt P_i^\pi \Bigg) + \mathtt C_k^\ell - \mathtt P_k
	\end{align}
	the net costs of all network participants in cell $k \in \mathcal Z$.
	Consider the cost minimization problems
	summarizing the optimization problems\footnote{Eq. \eqref{eq:Pareto-OP-price-constraint} contains all elements of \eqref{eq:extended-constraint-D+} that belong to the nodes in cell $k$.} of all network participants in cell $k$:
		\begin{subequations}\label{cell-cost-minimization}
		\begin{IEEEeqnarray}{Cll}
			\quad \min_{\bm p_{g,k}, \bm U_{f, k}, \bm U_{\mathcal I, k}}& \mathtt C_k^\mathrm{total}	\\
			\quad \text{s.t.} & (\bm D_c^{k+})^\top \bm \lambda_k = \nullvec{},\label{eq:Pareto-OP-price-constraint} \\
			& \underline{\bm p}_{g,k} \leq {\bm p}_{g,k} \leq \overline{\bm p}_{g,k}, \\
			 & \underline{\bm U}_{f,k} \leq {\bm U}_{f,k} \leq \overline{\bm U}_{f,k}, \\
			  & \underline{\bm U}_{\mathcal I,k} \leq {\bm U}_{\mathcal I,k} \leq \overline{\bm U}_{\mathcal I,k}. 
		\end{IEEEeqnarray}
	\end{subequations}
	Let $\bm \kappa' > \nullvec{}$ be fixed and the corresponding value of $\mathtt C_k^\mathrm{total}$ for an equilibrium of 
	\eqref{cell-cost-minimization}
	with parameter $\bm \kappa$ set to $\bm \kappa'$
be denoted by $\mathtt C_k^\mathrm{total \star}(\bm \kappa')$. Then, there exists no other \emph{dominating} $\bm \kappa''>\nullvec{}$ with $\bm \kappa'' \neq \bm \kappa'$ such that the following two conditions hold:
	\begin{subequations}\label{dominated-kappa}
	\begin{align}
	\forall \; k \in \mathcal Z: && \mathtt C_k^\mathrm{total \star}(\bm \kappa'') \leq \mathtt C_k^\mathrm{total \star}(\bm \kappa') \\ 
		\exists \; k \in \mathcal Z: && \mathtt C_k^\mathrm{total \star}(\bm \kappa'') < \mathtt C_k^\mathrm{total \star}(\bm \kappa')
	\end{align}
	\end{subequations}
	\end{proposition}
\begin{IEEEproof}
	Taking into account that \eqref{cell-cost-minimization} are convex optimization problems and following the same lines as in the proof of Proposition \ref{prop:equivalence-to-centralized-optimization}, it can be shown that each KKT point of \eqref{cell-cost-minimization} is given by \eqref{eq:central-closed-loop-begin}--\eqref{eq:central-closed-loop-end} and vice versa. 
	For each $\bm \kappa>\nullvec{}$, the equivalent problem \eqref{eq:central-OP} is a linear scalarization of \eqref{cell-cost-minimization} with positive weights $1 / \kappa_k >0$.
	With \cite[Proposition~9]{Emmerich.2018}, it is implied that for each $\bm \kappa'>\nullvec{}$, the solution of \eqref{cell-cost-minimization} is a part of the Pareto front,
	hence there exists no $\bm \kappa'>\nullvec{}$ which is dominated by another $\bm \kappa''>0$ in terms of \eqref{dominated-kappa}.
\end{IEEEproof}	

The remaining degrees of freedom manifest in a free choice of the participation factor $\bm \kappa \in \mathds R^{|\mathcal Z|}_{>0}$.
\begin{neu}{By adjusting $\bm \kappa$ appropriately, specific desired relationships between zonal prices are imposed,
which can serve as a real-time control mechanism for the reallocation of electricity supply from one cell to another.}\end{neu}

%
%
%
%
%
%
%
%
%
%
%
%
%
%
%
%
%

\begin{neu}
\section{Application: Flow-Based Congestion Management}\label{ch:application}
In the following, we exemplarily apply the automatic regulation of $\bm \kappa$ for a continuous-time congestion management strategy, which is conducted among neighboring CCs {by adjusting} the cell-specific prices.

For this, we describe the cell topology by the condensed graph ${\mathscr G}_z=(\mathcal Z,\widehat{\mathcal E})$, where the set $\widehat{\mathcal E}$ contains all \emph{inter-cell lines}, i.e. lines that connect nodes from different cells.
Let $P_{ij}$ denote the sending-end active power flow from node $i$ to $j$.
Since $P_{ij} \neq -P_{ji}$ for lossy lines, we define
\begin{align}
P_{m}=\begin{cases} P_{ij}, & |P_{ij}| \geq |P_{ji}|, \\
-P_{ji}, & \text{otherwise}
\end{cases}
\end{align}
for each line $m \in \widehat{\mathcal E}$ from $i$ to $j$.
With $P_m^{\max}$ specifying the maximum permissible active power flow of line $m \in \widehat{\mathcal E}$, the \emph{congestion rate} of line $m \in \widehat{\mathcal E}$ can thus be calculated as $\mathcal C_m = P_m / P^{\max}_m$. 
{Accordingly, a line is congested whenever $\mathcal C_m > 1$ or $\mathcal C_m < -1$.

To penalize lines which are close to congestion, we define a barrier function for each line as follows
\begin{align}
\gamma_m = \begin{cases}
\displaystyle\mathcal C_m \cdot 
\displaystyle\frac{|\mathcal C_m|-\mathcal C_m^{\min}}{(1-|\mathcal C_m|)(1-\mathcal C_m^{\min})}, & |\mathcal C_m| \geq \mathcal C_m^{\min}, \\ 
0, & |\mathcal C_m| < \mathcal C_m^{\min},
\end{cases}
\end{align}
where 
$0 < \mathcal C_m^{\min} < 1$ is a user-defined threshold for control actions. 
As illustrated in Fig. \ref{fig:barrierfunction}, $\gamma_m \to \pm \infty$ if $\mathcal C_m \to \pm 1$, and $\gamma_m =0$ if $|\mathcal C_m|$ is below its threshold $\mathcal C_m^{\min}$. 
\begin{figure}
		\setlength\abovecaptionskip{-0.2\baselineskip}
	\centering
	\includegraphics[width=0.925\columnwidth]{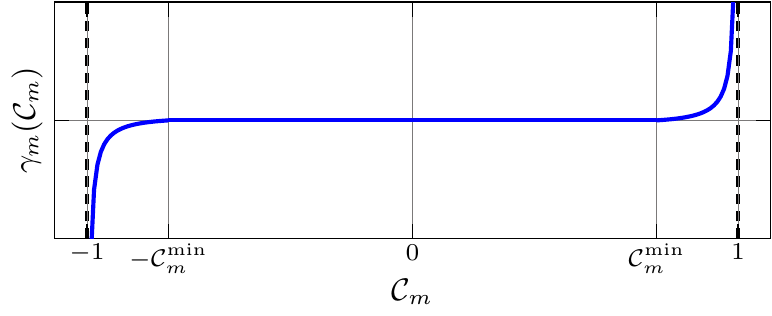}
	\caption{Plot of barrier function $\gamma_m(\mathcal C_m)$}
	\label{fig:barrierfunction}
\end{figure}

Firstly, consider a single congested line $m \in \widehat{\mathcal E}$ connecting two different cells $k_1$ and $k_2$. It is intuitive to counteract the congestion by stipulating a price difference between both cells. Hence, we set 
\begin{align}\label{cong-idee}
\ln\left(\frac{\kappa_{k_1}}{\kappa_{k_2}}\right) = -\gamma_m.
\end{align} 
This
leads to $\kappa_{k_1}<\kappa_{k_2} $, if $\mathcal C_m$ is above the threshold,
$\kappa_{k_1} > \kappa_{k_2}$ if $-\mathcal C_m$ is above the threshold, 
and $\kappa_{k_1}=\kappa_{k_2}$ otherwise. 

Next, we consider all inter-cell lines $m \in \widehat{\mathcal E}$.
With the shorthand notations} $\phi_k := \ln(\kappa_k)$, $\bm \phi = \mathrm{col}_k\{\phi_k\}$ and $\bm \gamma = \mathrm{col}_m\{\gamma_m\}$, \eqref{cong-idee} can be written in vector-matrix notation as
\begin{align}\label{cong-idee-matrix}
\bm \gamma = -{\bm D}_z^\top \bm \phi,
\end{align} 
where ${\bm D}_z$ is the incidence matrix of ${\mathscr G}_z$. 
Motivated by the fact that ${\mathscr G}_z$ is meshed and may have multiple (parallel) edges\footnote{\begin{neu}Note that each inter-cell line is represented by a separate edge in ${\mathscr G}_z$.\end{neu}}, the overall share of each cell on congestion of adjacent lines can be calculated by mapping the respective components in \eqref{cong-idee-matrix} to each cell, i.e. left-multiplying with ${\bm D}_z$, which yields
\begin{align}\label{eq:phi-kappa}
{\bm D}_z\bm \gamma = -\underbrace{{\bm D}_z {\bm D}_z^\top}_{=:\mathcal B}\bm \phi,
\end{align}
where $\mathcal B$ is the Laplacian matrix of ${\mathscr G}_z$.
{To calculate $\bm \phi$ in \eqref{eq:phi-kappa} as well as the resulting} participation factors $\bm \kappa$, we adapt
\eqref{eq:phi-kappa} by a gradient descent, which finally leads to { the controller equations} 
\begin{subequations}\label{eq:kappa-congestion-regler}
\begin{align}
\bm \tau_{\bm \phi} \dot{\bm \phi} &= - {\mathcal B}\bm \phi - {\bm D}_z\bm \gamma, \label{phi-congestion-reglergleichung}\\
\bm \kappa &= \mathbf{exp}\{\bm \phi\}
\end{align}
\end{subequations}
with $\bm \tau_{\bm \phi} > \nullvec{}$.
\begin{remark}
Note that the individual values of $\phi_k$ and $\kappa_k$ in \eqref{eq:kappa-congestion-regler} can be calculated locally by the respective CC, only requiring information about neighboring $\kappa$ values and adjacent power flows across its own cell boundary.
\end{remark}
\end{neu}


\section{Case Study}\label{ch:woc-simulation}
\begin{figure}
	\centering
	\includegraphics[width=\columnwidth]{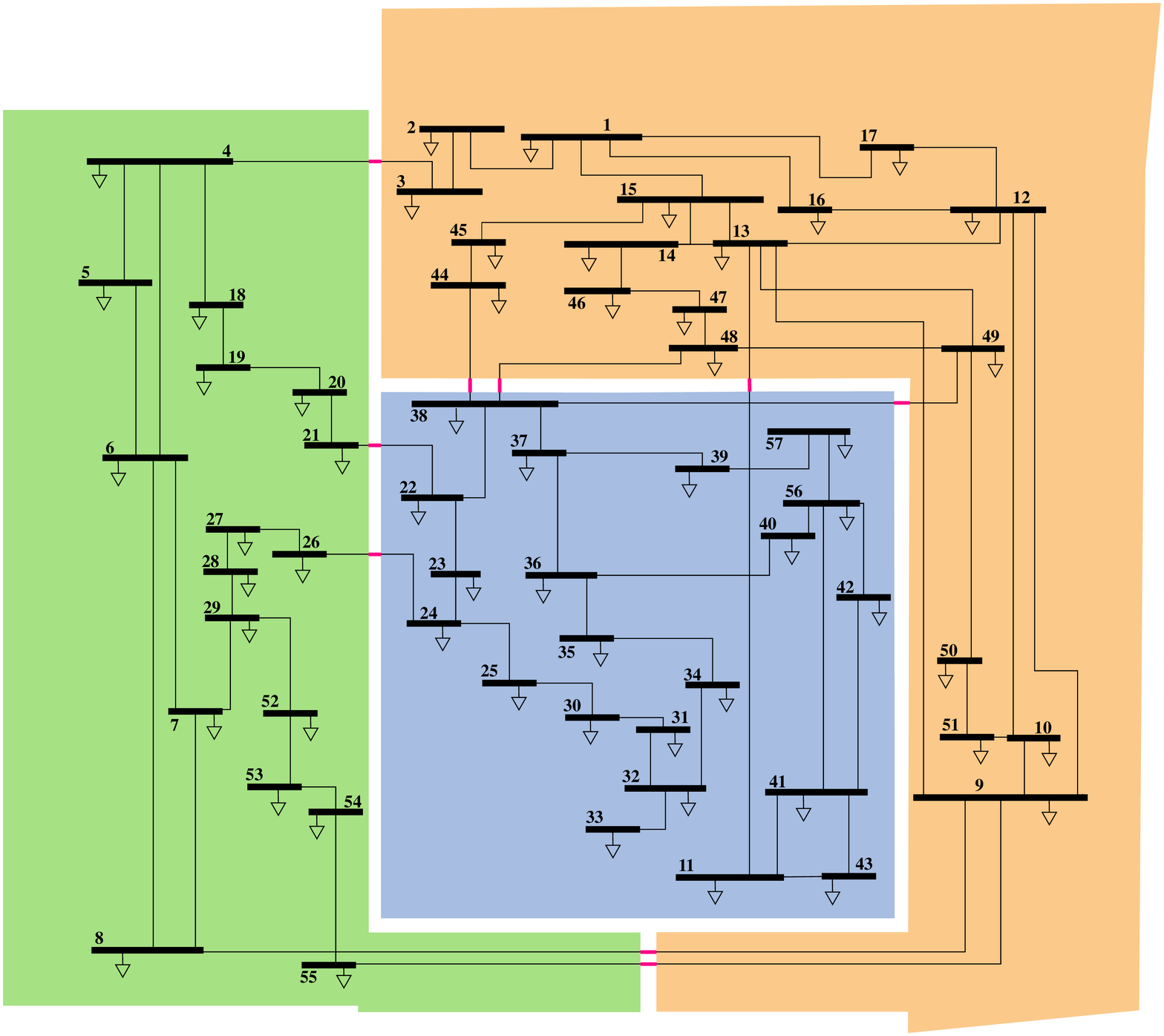}
	\caption{IEEE 57-bus system divided into \begin{neu}three\end{neu} cells \cite{IEEE57.Online}.}
		\label{fig:IEEE57}
	\end{figure}
We validate the performance of the developed distributed controller on the IEEE 57-bus test case, which is divided into three cells (see Fig. \ref{fig:IEEE57}).
The assignment \begin{neu}of buses\end{neu} to the node types is \begin{neu}given in\end{neu} Table \ref{tab:nodetypes-IEEE-case}.
\begin{neu}\subsection{Scenarios}\end{neu}
Let $\bm D_{\text{IEEE}57}$ denote the incidence matrix of the IEEE 57-bus test case
and $\bm D_{\text{IEEE}57}^{-}$ denote the respective incidence matrix if all boundary edges of adjacent cells \begin{neu}(colored magenta in Fig. \ref{fig:IEEE57})\end{neu} are removed.
To investigate the role of connectivity of $\mathscr G_p$ and $\mathscr G_c$ and the effects of $\bm \kappa$, we examine four different topological scenarios:
\begin{neu}
\begin{enumerate}[(I)]
	\item \emph{Isolated cells:} All physical connections via inter-cell lines are removed and the communication networks of CCs are separated from each other, i.e. $\bm D_p=\bm D_c=\bm D_{\text{IEEE}57}^-$.
	\item \emph{Free zonal prices:} The cells are physically connected, but the communication networks of CCs are separated from each other, i.e. $\bm D_p=\bm D_{\text{IEEE}57}$ and $\bm D_c=\bm D_{\text{IEEE}57}^-$. 
	\item \emph{Uniform prices:} Both the physical lines and the communication networks of cells are connected and the participation factors are set to a constant value of one, i.e. 
	$\bm D_p = \bm D_c = \bm D_{\text{IEEE}57}$ with $\bm \kappa = \mathbb 1$.
	\item \emph{Uniform prices + congestion management:} Both the physical lines and the communication networks of cells are connected, i.e.
	$\bm D_p = \bm D_c = \bm D_{\text{IEEE}57}$, and the participation factors are set according to the real-time congestion management strategy \eqref{eq:kappa-congestion-regler}.
\end{enumerate}
\end{neu}

\subsection{Network Model and Parameterization}
All of the following numerical values are given in p.u. with $U_{\mathrm{base}}=\SI{135}{\kilo\volt}$,
$S_{\mathrm{base}}=\SI{100}{{\mega\volt\ampere}}$,
\begin{neu}and $\mathtt{C}_{\mathrm{base}}= 1\,\text{MU} / S_{\mathrm{base}}$.
The controller parameters $\bm \tau_g$ in \eqref{eq:WOC-closed-loop-begin} are set to \num{0.1}, while the controller parameters $\bm \tau_{(\cdot)}$ in \eqref{eq:WOC-closed-loop-zweites}--\eqref{eq:WOC-closed-loop-end} are set to \num{0.01}.\end{neu}
\begin{table}[!t]
	\renewcommand{\arraystretch}{1.1}
	\caption{Assignment of Node Types}
	\centering
	\begin{tabular}{|c||l|}
		\hline
		$\mathcal G$ & 2,\,3,\,6,\,8,\,9,\,10,\,19,\,21,\,29,\,30,\,32,\,34,\,37,\,39,\,40,\,41,\,44,\,48,\,55\\
		$\mathcal I$ & 4,\,11,\,14,\,15,\,16,\,17,\,18,\,22,\,24,\,25,\,26,\,33,\,36,\,42,\,45,\,46,\,49,\,50,\,53\\
		$\mathcal L$ &1,\,5,\,7,\,12,\,13,\,20,\,23,\,27,\,28,\,31,\,35,\,38,\,43,\,47,\,51,\,52,\,54,\,56,\,57\\
		\hline
	\end{tabular}
\label{tab:nodetypes-IEEE-case}
\end{table}
The node parameters are chosen to be randomly distributed within the specific intervals in Table \ref{tab:node-parameters}.
\begin{table}[!t]
	\renewcommand{\arraystretch}{1.1}
	\caption{Node Parameters}
	\centering
	\begin{tabular}{|c||c|c|c|c|}
		\cline{2-5}
		\multicolumn{1}{c|}{}& $A_i$ & $M_i$ & $X_{d,i} - X_{d,i}'$ & $\tau_{U,i}$ \\
		\hline
		$\mathcal G$ & $[1.2; 1.7]$ & $[20;27]$ & $[0.12;0.19]$ & $[6.4;7.7]$ \\
		$\mathcal I$ & $[1.2;1.7]$ & $[4;5.5]$ & --- & --- \\
		$\mathcal L$ & $[1.2;1.7]$ & --- & --- & --- \\
		\hline
		\end{tabular}
	\label{tab:node-parameters}
\end{table}
The lower and upper bounds of $\bm p_g$ are set to $\underline p_{g,i}=-0.002$ and $\overline p_{g,i}=0.003$ respectively, and the voltage limits are set to $\underline U_i = 0.98$ and $\overline U_i = 1.02$. 
\begin{neu}For scenario IV, we specify a maximum power flow of  $P_m^{\max}= \num{0,01}$ and a threshold of $\mathcal C_m^{\min}=\num{0,8}$ for 
all inter-cell lines $m \in \widehat{\mathcal E}$ and set the controller parameters $\bm \tau_\phi$ in \eqref{phi-congestion-reglergleichung} to \num{10}. \end{neu}

At $t=0$, the system is in synchronous mode with $\bm L=\nullvec{}$. All controller variables $\bm p_g$, $\bm \nu$, $\bm \lambda$, $\bm \mu_{(\cdot)}$ are initialized to zero. 
Without loss of generality, we set $\mathcal P \equiv \mathcal Z$ and choose the cost functions to
\begin{align}
\mathtt C_{\pi}(\bm p_\pi) = \frac 12 \cdot \bigg(\sum_{i \in \mathcal V_\pi^{\mathcal P}} \frac{1}{w_i} \cdot p_{g,i}^2 \bigg), && \pi \in \mathcal P,\label{IEEE-57-cost-function}
\end{align} 
where $w_i = 1 + 0.04 \cdot (i-1)$.
The initial values of 
$\vartheta_{ij}$ are chosen within the interval $[-0.04;0.014]$ and all voltages $U_i$ are chosen within the interval $[0.98;1.02]$.

The active and reactive power consumptions $\bm p_\ell$ and $\bm q_\ell$ are piecewise constant in order to simulate a step change in generation or consumption in specific areas of the network. To this end, the initial values of $\bm p_\ell$ and $\bm q_\ell$  are chosen such that the AC power flow equations are satisfied.
Subsequently, the following stepwise load \begin{neu}changes\end{neu} are applied:
\begin{itemize}
	\item at $t_1=\begin{neu}\SI{5}{\minute}\end{neu}$, the active power consumption $p_{\ell,28}$ in cell 3 increases by \begin{neu}$\num{0.015}$\end{neu}, 
	\item at $t_2=\begin{neu}\SI{10}{\minute}\end{neu}$, the reactive power consumption $q_{\ell,28}$ at the same node increases by \begin{neu}$\num{0.015}$\end{neu}, 
	\item at $t_3=\begin{neu}\SI{15}{\minute}\end{neu}$, two additional nodes, 20 and 27, which are also located in cell 3, increase their active and reactive power consumption by \begin{neu}\num{0,0075}\end{neu},
	\item at $t_4=\begin{neu}\SI{20}{\minute}\end{neu}$, the power consumptions of nodes 20, 27, and 28 are reset to their initial values. At the same time, active and reactive power consumption at load nodes  12, 13, and 43 in cells 1 and 2 increases by \begin{neu}$\num{0,0075}$\end{neu} each,
	\item at $t_5 =\begin{neu}\SI{25}{\minute}\end{neu}$ the active and reactive power consumptions of all above-mentioned nodes are multiplied by $-1$ to simulate a reversed load flow.
\end{itemize}
\begin{neu}
The numerical implementation is conducted on a machine with an Intel Core i7-6600U and 12 GB of RAM using the numerical solver \texttt{NDSolve} in Wolfram Mathematica (Version 12.0.0).
\subsection{Numerical Results}
Figs. \ref{fig:scenario-B3}--\ref{fig:scenario-Kappa1Cong} show the resulting nodal prices $\bm \lambda$, active power generations $\bm p_g$, nodal frequencies $\bm f$ and voltage magnitudes $\bm U$ for scenarios I--IV, respectively.
All four scenarios share the similarity that no 
price differences can be observed within each cell $k \in \{1,2,3\}$, see the upper plot in Figs. \ref{fig:scenario-B3}--\ref{fig:scenario-Kappa1Cong}. This shows that the nodal prices have a very fast convergence to a cell-specific price $\Lambda_k$.
\begin{figure}
	\begin{flushright}
		\includegraphics[width=0.98\columnwidth]{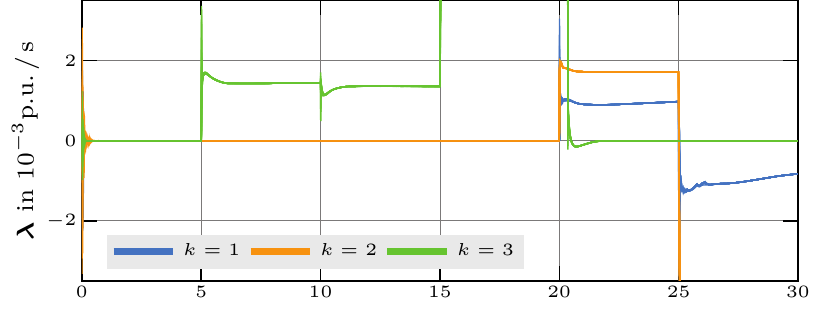}
		\includegraphics[width=0.985\columnwidth]{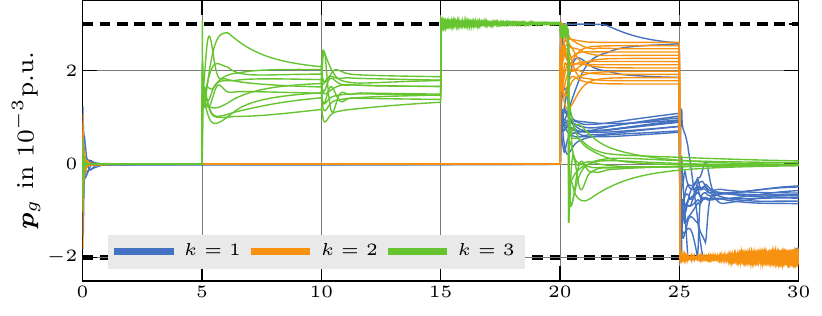}
		\includegraphics[width=\columnwidth]{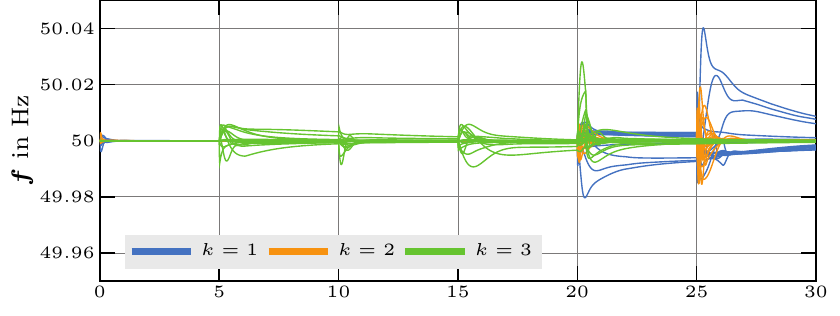}
		\includegraphics[width=0.985\columnwidth]{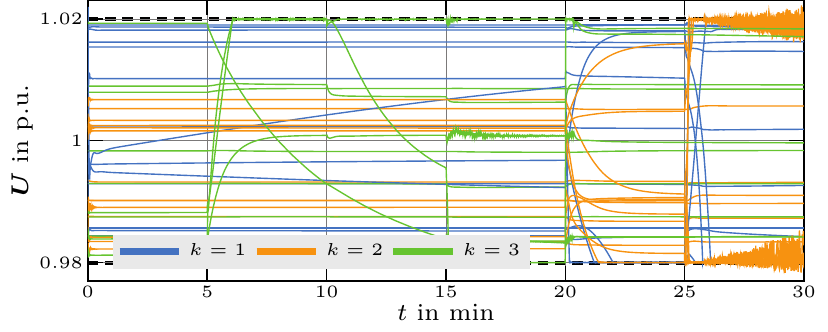}
	\end{flushright}
	\caption{\begin{neu}Simulation results for scenario I (isolated cells)\end{neu}.}
	\label{fig:scenario-B3}
\end{figure}
\subsubsection*{Scenario I (isolated cells)}
It can be seen in the first subplot of Fig. \ref{fig:scenario-B3}, that
due to physical decoupling,
the loads steps in cell 3 at $t_1$, $t_2$, and $t_3$ only affect the price in cell 3, while
prices in other cells remain unaffected.
Likewise, the load step in cells 1 and 2 at $t_5$ affects the prices in cells 1 and 2, while the price in cell 3 remains to its initial value.
With regard to the active power injections, it can be seen that the upper and lower limits 
$\overline{\bm p}_g$ and $\underline{\bm p}_g$ are met except after $t_5$ in cell 2, where we see a significant oscillation of all injections around the lower limit.
Moreover, the active power injections in cell 1 do not converge within the interval between $t_4$ and $t_5$.
As evident in the third subplot of Fig. \ref{fig:scenario-B3},
each load step causes temporary deviation of the nodal frequencies in the involved cells, which are subsequently regulated back to their nominal value. 
The maximum overshoot across the entire simulated time is $+\SI{0.04}{\hertz}$.
After the steps at $t_4$ and $t_5$, the performance of the frequency regulation in cell 1 is poor. 
The last subplot of Fig. \ref{fig:scenario-B3} shows the resulting voltage magnitudes in all cells. Again, all upper and lower limits are met except after $t_5$, where most of the voltages in cell 2 oscillate around its upper or lower limits.
{
This indicates that in the isolated scenario, PPOs in cell 2 are unable to meet the local active and reactive power demand without violating the constraints \eqref{eq:PPO-OP-constraint-begin}--\eqref{eq:PPO-OP-constraint-end}.	
}

\begin{figure}
	\begin{flushright}
		\includegraphics[width=0.98\columnwidth]{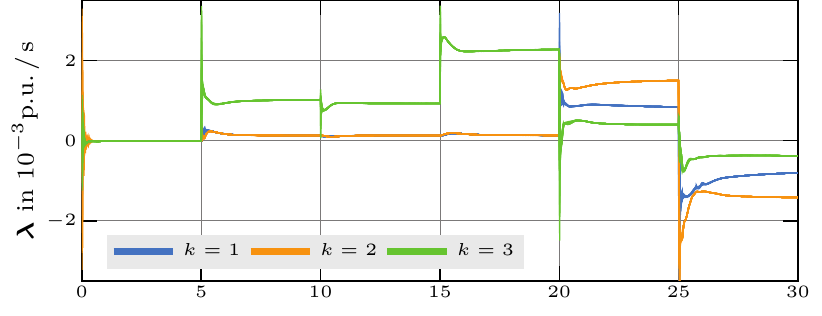}
		\includegraphics[width=0.985\columnwidth]{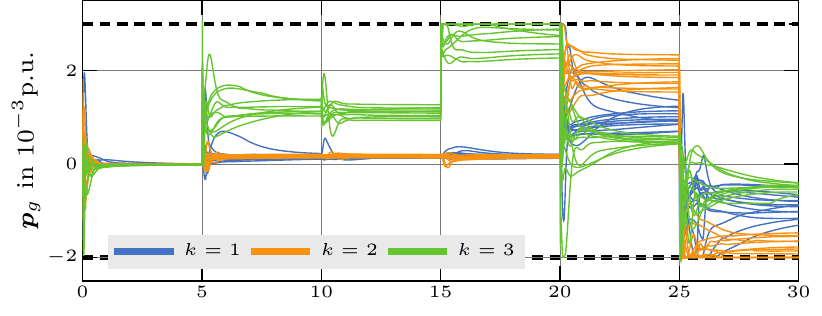}
		\includegraphics[width=\columnwidth]{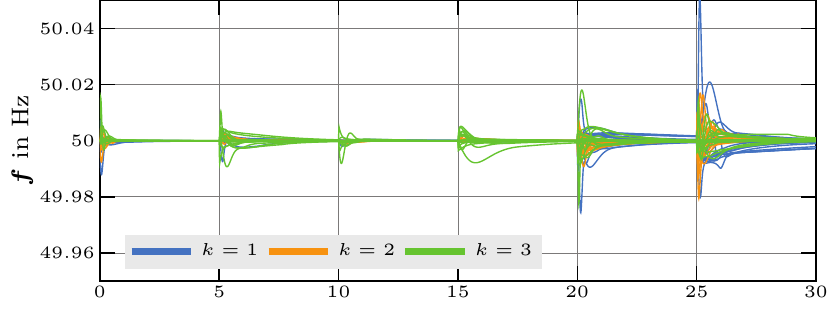}
		\includegraphics[width=0.985\columnwidth]{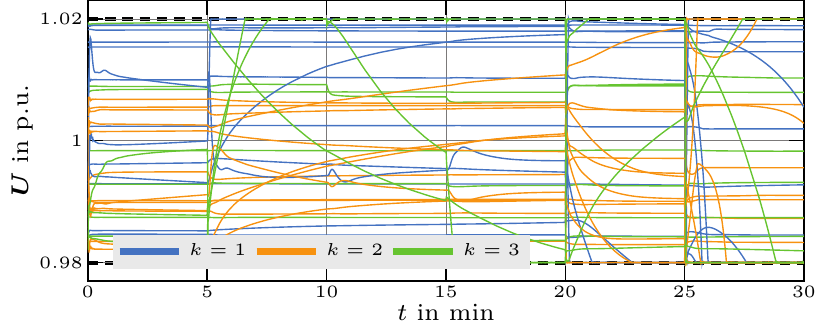}
	\end{flushright}
	\caption{\begin{neu}Simulation results for scenario II (free zonal prices)\end{neu}.}
	\label{fig:scenario-A3}
\end{figure}
\subsubsection*{Scenario II (free zonal prices)}
The resulting cell-specific prices shown in the upper plot in Fig. \ref{fig:scenario-A3} 
predominantly follow a similar course as scenario I. In particular, the load steps in cell 3 at $t_1$, $t_2$, and $t_3$ lead to an increase in the price in cell 3. Likewise, the load steps in cells 1 and 2 at $t_4$ lead to an increase in the price in cells 1 and 2.
However, in contrast to scenario I, load steps in one cell also affect the prices in the other cells, signifying the effect of physical interconnection via inter-cell lines.
Moreover, the spread between zonal prices is lower than in scenario I.
Similar to the price curves, the resulting active power injections in different cells are affected by each other and show a lower spread than in scenario I.
This illustrates how the other cells contribute to eliminating the imbalance in a specific cell. 
The nodal frequencies reveal a maximal overshoot of $+\SI{0.05}{\hertz}$ but better convergence as compared to scenario I.
Lastly, voltage limits, just like all other limits, are maintained without oscillation.

\begin{figure}
	\begin{flushright}
		\includegraphics[width=0.98\columnwidth]{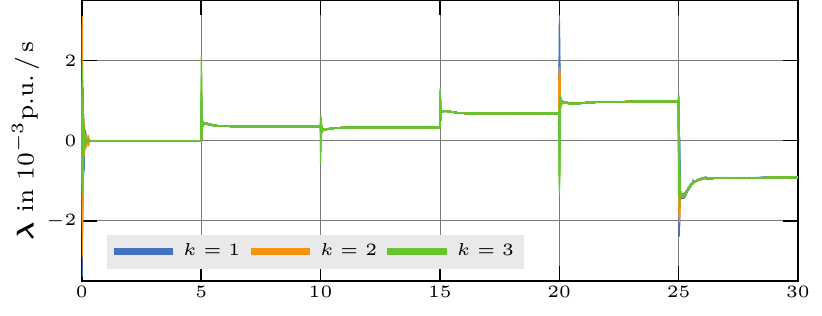}
		\includegraphics[width=0.985\columnwidth]{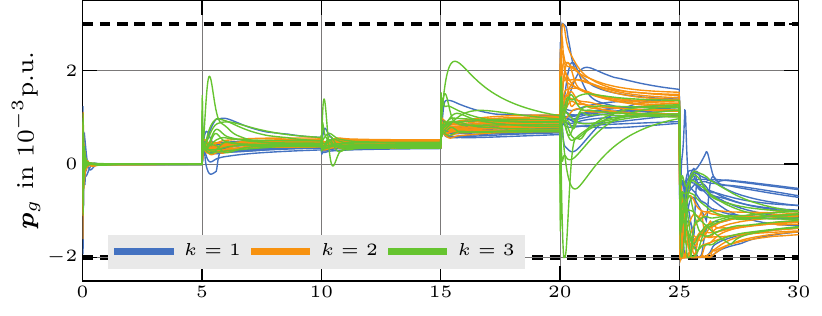}
		\includegraphics[width=\columnwidth]{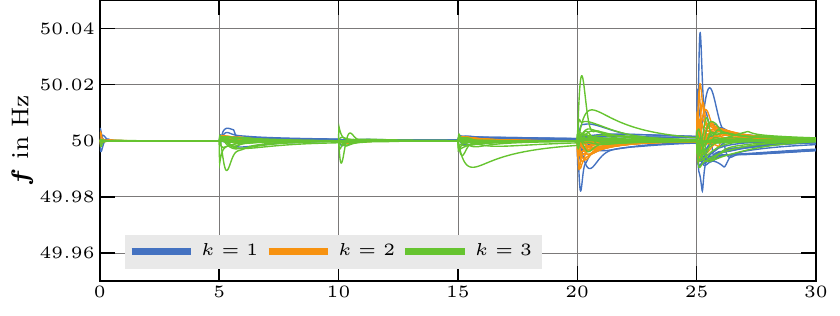}
		\includegraphics[width=0.985\columnwidth]{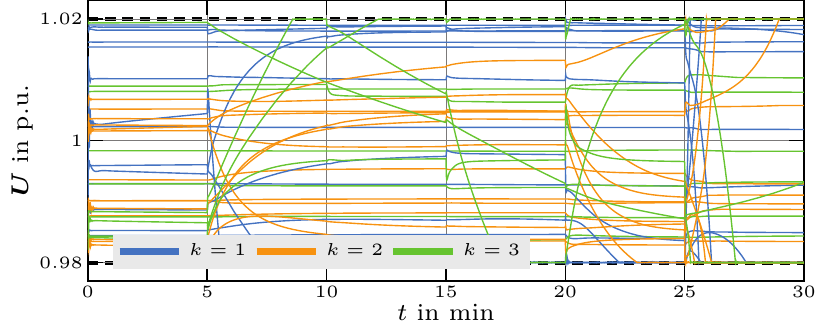}
	\end{flushright}
	\caption{\begin{neu}Simulation results for scenario III (uniform prices)\end{neu}.}
	\label{fig:scenario-Kappa1}
\end{figure}
\subsubsection*{Scenario III (uniform pricing)}
The upper plot in Fig. \ref{fig:scenario-Kappa1} reveals that in this scenario, all prices quickly synchronize to a common value after each load step.
Accordingly, load steps in one specific cell impact the active power injections, nodal frequencies, and voltage magnitues in all cells.
This corresponds to the individual active power injections being closer to each other, barely reaching its specific upper or lower limits.
The resulting nodal frequencies are similar to those of scenario II.
Again, voltage magnitudes as well as all other variables are kept within their respective limits.

\begin{figure}
	\begin{flushright}
		\includegraphics[width=0.98\columnwidth]{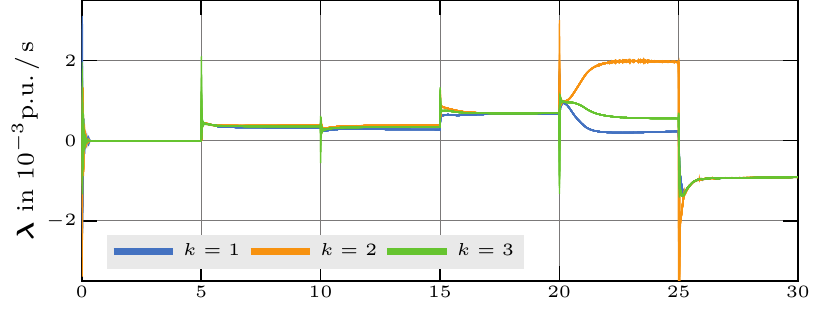}
		\includegraphics[width=0.985\columnwidth]{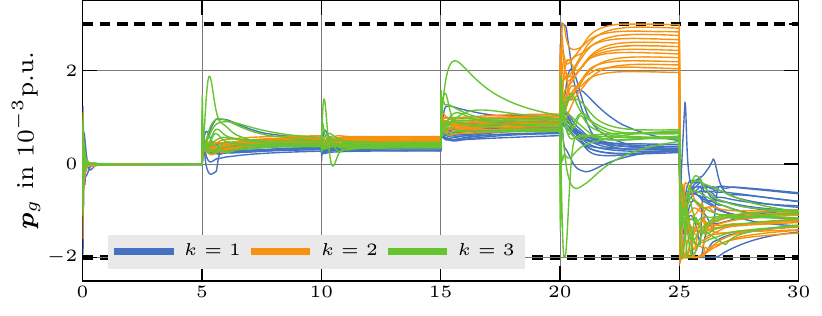}
		\includegraphics[width=\columnwidth]{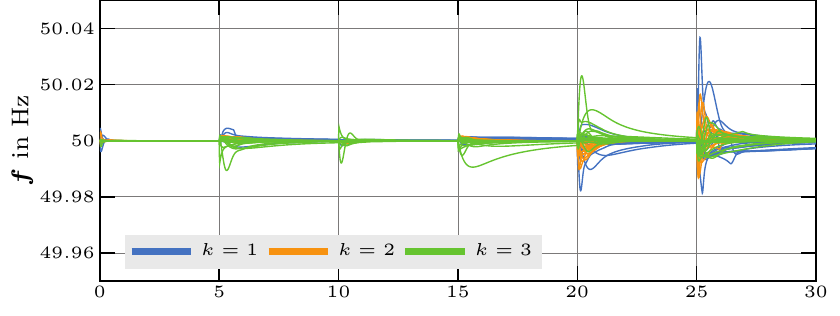}
		\includegraphics[width=0.985\columnwidth]{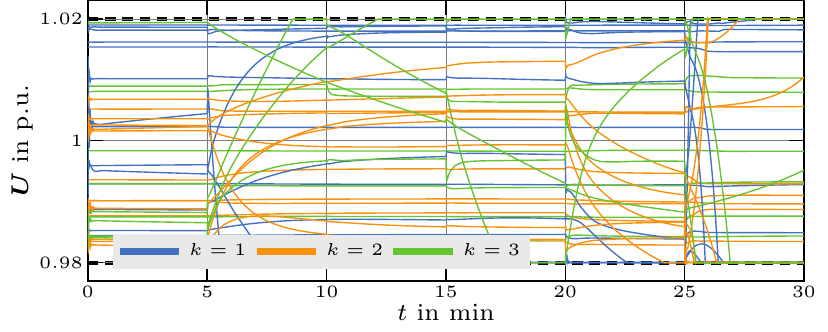}
	\end{flushright}
	\caption{\begin{neu}Simulation results for scenario IV (uniform pricing + congestion management)\end{neu}.}
	\label{fig:scenario-Kappa1Cong}
\end{figure}

\begin{figure}
	\centering
	\includegraphics[width=\columnwidth]{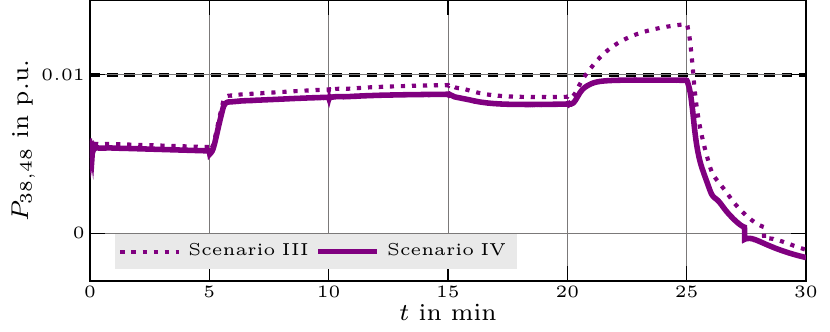}
	\caption{\begin{neu}Resulting power flow on line $(38,48)$ for scenarios III and IV\end{neu}.}
	\label{fig:congline}
\end{figure}

\subsubsection*{Scenario IV (uniform pricing + congestion management)}
As can be seen in Fig. \ref{fig:scenario-Kappa1Cong}, the resulting cell-specific prices are similar to those of scenario III. However, after the load step at $t_4$, there is a considerable split into three different zonal prices. Subsequently, after $t_5$, prices re-synchronize to a common value, which is equal to that of scenario III.
This result can be made plausible as follows:
In scenario III, all resulting inter-cell power flows are below its specified maximum $P_m^{\max}$, except for line $(38,48)$ between cells 1 and 2 which is overloaded (see dotted line in Fig. \ref{fig:congline}).
Accordingly, in scenario IV the real-time congestion management becomes active, which leads to a non-uniform $\bm \kappa(t)$ between $t_4$ and $t_5$. In particular, $\kappa_2$ is increased significantly, while $\kappa_1$ decreases significantly and $\kappa_3$ decreases slightly. The overall active power injection is thus shifted from cell 1 to cell 2, which results in $P_{38,48}$ falling right below its maximum specified value (see solid line in Fig. \ref{fig:congline}). After $t_5$, the congestion management becomes inactive, hence prices automatically synchronize
to the same common value as in scenario III.
Concurrently, nodal frequencies and voltages are not compromised by the splitting of prices.

\subsection{Computational Performance}
Table \ref{tab:computingtime} shows the resulting simulation times for the numerical solver,
ranging from \SI{67,86}{\second} for scenario III up to \SI{278,05}{\second} for scenario I.
All values are well below the simulated time of \SI{1800}{\second}, indicating the real-time applicability of the proposed method under each scenario.
Table \ref{tab:computingtime} as well as additional studies suggest that the computation time of the solver increases with the number of variables that reach their limits.
However, since the proposed control framework is distributed, in real-world applications, separate lower-order controllers will be operated in parallel at every single node, thereby resulting in significantly lower computation times.
\end{neu}

\begin{table}[!t]
	\renewcommand{\arraystretch}{1.1}
	\caption{\begin{neu}Comparison of Simulation Times\end{neu}}
	\centering
	\begin{tabular}{|l||c|c|c|c|}
		\hline
		Scenario & I & II & III & IV \\
		\hline
		Simulation time in \si{\second} & \num{278.65} & \num{102,40} & \num{67.86} & \num{166.05}\\
		\hline
	\end{tabular}
	\label{tab:computingtime}
\end{table}

\section{Conclusion}\label{ch:woc-conclusion}
In this paper, we presented a distributed control strategy for 
\begin{neu}real-time\end{neu} 
dynamic pricing in zonal electricity markets.
Incorporating the WoC concept, the interplay of different optimization problems of PPOs and CCs resulted in an overall control system for which frequency and voltage \begin{neu} regulation is always preserved.\end{neu} 
Simulation studies on a modified IEEE 57-bus system compared an \begin{neu}isolated\end{neu} operation with different zonal pricing approaches.
\begin{neu}
	While in the first scenario, local generation capacity was shown to be unable to meet the local demand 
	in all cells (scenario I), the interconnection of price zones led to an
	improved participation of neighboring PPOs from adjacent cells, which resulted in a lower price spread (scenario II).
	Through the introduction of the participation factor $\bm \kappa$, specific relationships between zonal prices can be imposed which may either be constant in time to achieve uniform prices throughout the network (scenario III),
	or time-varying to  
	relieve congestion at market clearing with an automatic settlement of regionally differentiated prices in case of heavily loaded lines {(scenario IV)}.\end{neu} 
The \begin{neu}overall\end{neu} controller features low efforts for parameter tuning, since the choice of the free controller parameters $\bm \tau_{(\cdot)}$ 
does not affect the value of the closed-loop equilibrium.


\begin{neu}By permitting automatic control of $\bm \kappa$, our work has given rise to allocating marginal production in specific regions of the power system in real time.
Further extensions to the presented congestion management strategy may be undertaken by  
including
{alternative short- and long-term goals into the generation of real-time prices.
These findings could help 
in providing a new perspective on the optimal partitioning of price zones.}
\end{neu}

%
%
%
%
%


%
%
%
%

\todo[disable]{
\begin{underconstruction}
\section{Resterampe}
Through regionally differentiated prices, the grid-awareness of marginal production in a specific location can be valuated

Congestion is avoided at market clearing

		spikes in generation caused by intermittent renewable generation that could lead to violations of the capacity constraints.

The cell-specific prices are settled in real-time by means of neighbor-to-neighbor communication between CCs.

\end{underconstruction}
}


%
%
%
%
%
%
%
%
%
%


%

\appendices


\ifCLASSOPTIONcaptionsoff
  \newpage
\fi



\bibliographystyle{IEEEtran}
\bibliography{bib/bib}
\end{document}